\documentclass[11pt]{amsart}
\usepackage{mathtools}
\usepackage{amssymb}
\usepackage[margin=1in]{geometry}
\usepackage[utf8]{inputenc}
\usepackage{mathtools}
\usepackage[rgb,usenames,dvipsnames]{xcolor}
\colorlet{linkcolour}{blue!50!black}

\RequirePackage{hyperref}
\hypersetup{
hypertexnames=false,
colorlinks=true,
linkcolor=linkcolour,
citecolor=linkcolour,
urlcolor=linkcolour,
}

\mathtoolsset{showonlyrefs}

\usepackage[ruled]{algorithm}
\usepackage[noend]{algpseudocode}

\usepackage{float} 
\usepackage{graphicx}
\usepackage{epstopdf}
\usepackage{color}
\usepackage{amsthm}
\theoremstyle{plain}
\newtheorem{theorem}{Theorem}
\newtheorem{corollary}[theorem]{Corollary}
\newtheorem{lemma}[theorem]{Lemma}
\newtheorem{proposition}[theorem]{Proposition}
\newtheorem{definition}[theorem]{Definition}

\newtheorem{remark}[theorem]{Remark}

\newcommand{\tangent}{\mathbf{t}}
\newcommand{\normal}{\mathbf{n}}
\newcommand{\binormal}{\mathbf{b}}
\newcommand{\darboux}{\mathbf{D}}
\newcommand{\dorbaux}{\mathbf{E}}
\newcommand{\spacey}{\mathcal{X}}
\newcommand{\opey}{\mathcal{U}}
\newcommand{\flag}{\mathbf{r}}
\newcommand{\flagsy}{\mathcal{F}}
\newcommand{\flagsyodd}{\flagsy_o}
\newcommand{\varpar}{\zeta}
\newcommand{\varfield}{W}
\newcommand{\llangle}{\langle\!\langle}
\newcommand{\rrangle}{\rangle\!\rangle}

\begin{document}
\title{Isometric Immersions and the Waving of Flags}
\author{Martin Bauer}
\address{Department of Mathematics, Florida State University}
\email{bauer@math.fsu.edu}
\author{Jakob M\o{}ller-Andersen}
\address{Department of Mathematics, Florida State University}
\email{jmoeller@math.fsu.edu}
\author{Stephen C. Preston}
\address{Department of Mathematics, Brooklyn College of City University New York}
\email{stephen.preston@brooklyn.cuny.edu}

\begin{abstract}
In this article we propose a novel geometric model to study the motion of a physical flag. In our approach a flag is viewed as an isometric immersion from the square with values in $\mathbb R^3$ satisfying certain boundary conditions at the flag pole. Under additional regularity constraints we show that the space of all such flags carries the structure of an infinite dimensional manifold and can be viewed as a submanifold of the space of all immersions. In the second part of the article we equip the space of isometric immersions with its natural  kinetic energy and derive the corresponding equations of motion.  This approach can be viewed in a similar spirit as
Arnold's geometric picture for the motion of an incompressible fluid.
\end{abstract}
\maketitle
\section{Introduction}

In this article we propose a geometric framework to model the motion of physical flags.
Mathematically, a flag on a flagpole may be modeled as an isometric $C^2$-immersion of a square into $\mathbb{R}^3$ subject to the constraint that one edge is mapped to the pole.
To obtain the simplest possible model, we ignore external forces and model the flag as though it follows a geodesic in the space of isometric immersions with Riemannian metric determined by the physical kinetic energy. Although the local problem of deformability of isometric immersions is well-known, see e.g., \cite{Spi1979} and the references therein, an additional difficulty in our setup comes from our need for \emph{global} coordinates, as well as the boundary conditions: matching the flag to the pole on one hand, and describing the edges of the square on the other hand. 
See Figure~\ref{fig:flagexamples} for three examples of flags (isometric immersions).

\subsection{Modelling equations of motions as geodesic equations}
Our approach follows similar geometrical models for other situations, such as modeling the motion of ideal fluids as a geodesic evolution in the group of volume-preserving diffeomorphisms, as first done by Arnold~\cite{Arn2014} in 1966. The advantage of this formulation is that it allows us to relate curvature of the manifold to stability of the system, and that it reduces the system to the simplest set of assumptions (without incorporating the details of external forces or the physical composition of the system). Another advantage is that it can lead to rigorous proofs of existence and uniqueness theorems by turning a PDE into an ODE, as in Ebin-Marsden~\cite{EM1970} for the incompressible Euler-equation. Many other PDEs have been recast as geodesics in various spaces, especially on diffeomorphism groups.
See \cite{KLMP2013,Kol2017,TiVi2011} for survey articles on the topic and Arnold-Khesin~\cite{AK1998} for an introduction to the field and more examples.

Diffeomorphism groups arise in studying motion of fluids which fill up their domain. In many other cases the system is a material moving in a higher-dimensional space, which leads one to work with spaces of embeddings and immersions; see \cite{BBM2014} and the references therein.  
{One example is the} motion of inextensible threads in Euclidean space: assuming the geometric constraint that the curve $\boldsymbol{\eta}$ preserves arc length $s$, we have the equation
\begin{equation}
\frac{\partial^2 \boldsymbol{\eta}}{\partial t^2} = \frac{\partial}{\partial s}\Big( \sigma \, \frac{\partial \boldsymbol{\eta}}{\partial s}\Big), \qquad \frac{\partial^2 \sigma}{\partial s^2} - \Big\lvert \frac{\partial^2 \boldsymbol{\eta}}{\partial s^2}\Big\rvert^2 \sigma = -\Big\lvert \frac{\partial^2 \boldsymbol{\eta}}{\partial t\partial s}\Big\rvert^2, \qquad \left\lvert \frac{\partial \boldsymbol{\eta}}{\partial s}\right\rvert^2 \equiv 1.
\end{equation}
{Here the function $\sigma$, which is determined by the second order ODE in the middle, can be interpreted as the tension of the curve. }  

This is a very old equation, but the existence theory is rather recent, as is the geometric treatment {(see ~\cite{preston2011motion,preston2012geometry} and references therein, along with the more recent \cite{csengul2017generalized} for weak solutions)}. A natural extension of this to higher dimensions is to consider the space of embeddings of surfaces into $\mathbb{R}^3$ with some constraint: either preserving the area element or preserving the Riemannian metric. Those that preserve the area element, which serve as a model for the motion of membranes in biological systems, were studied by several researchers including the first author~\cite{bauer2016riemannian,GaVi2014,molitor2017remarks}. In this article we study those that preserve the metric, which can serve as a model for unstretchable fabric or paper.

\subsection{Contributions of the article}
{For a general} {surface} {$M$ the study of the space of isometric immersions $M\to \mathbb R^3$ comes with several difficulties. The main reason for this is that }this space is relatively small and depends delicately on the geometry of the surface. For example, the Cohn-Vossen theorem says that if a {closed} surface has nonnegative Gaussian curvature and no open set where the the curvature is zero, then it is rigid: the only deformations are isometries of $\mathbb{R}^3$. For a recent survey of such results, see Han-Hong~\cite{han2006isometric}. If the Gaussian curvature is zero, as in our case, then there is a nontrivial family of deformations, but it is not very large: in the space of all immersions (three functions of two variables), the isometric immersions are generically described by two functions of one variable. Even this result is only valid locally, and we derive our own version. {Furthermore,} to work with the geodesic equation in this space, one  {would like to have}  a manifold structure, which cannot be expected for {general} {surfaces} { $M$}: in \cite{Weg1990} it has been shown that the space of isometric immersions { from $S^2$ to $\mathbb R^3$ is} not locally arcwise connected and hence cannot be a manifold, which gives a counter example to an earlier result of~\cite{Bin1984}. For the situation of this article---{isometric immersions from the flat square into $\mathbb R^3$}---we are nevertheless able to overcome these difficulties and show a manifold result, as described below.  
\begin{figure}[tbp]
	\begin{center}
		\includegraphics[width=.3\textwidth]{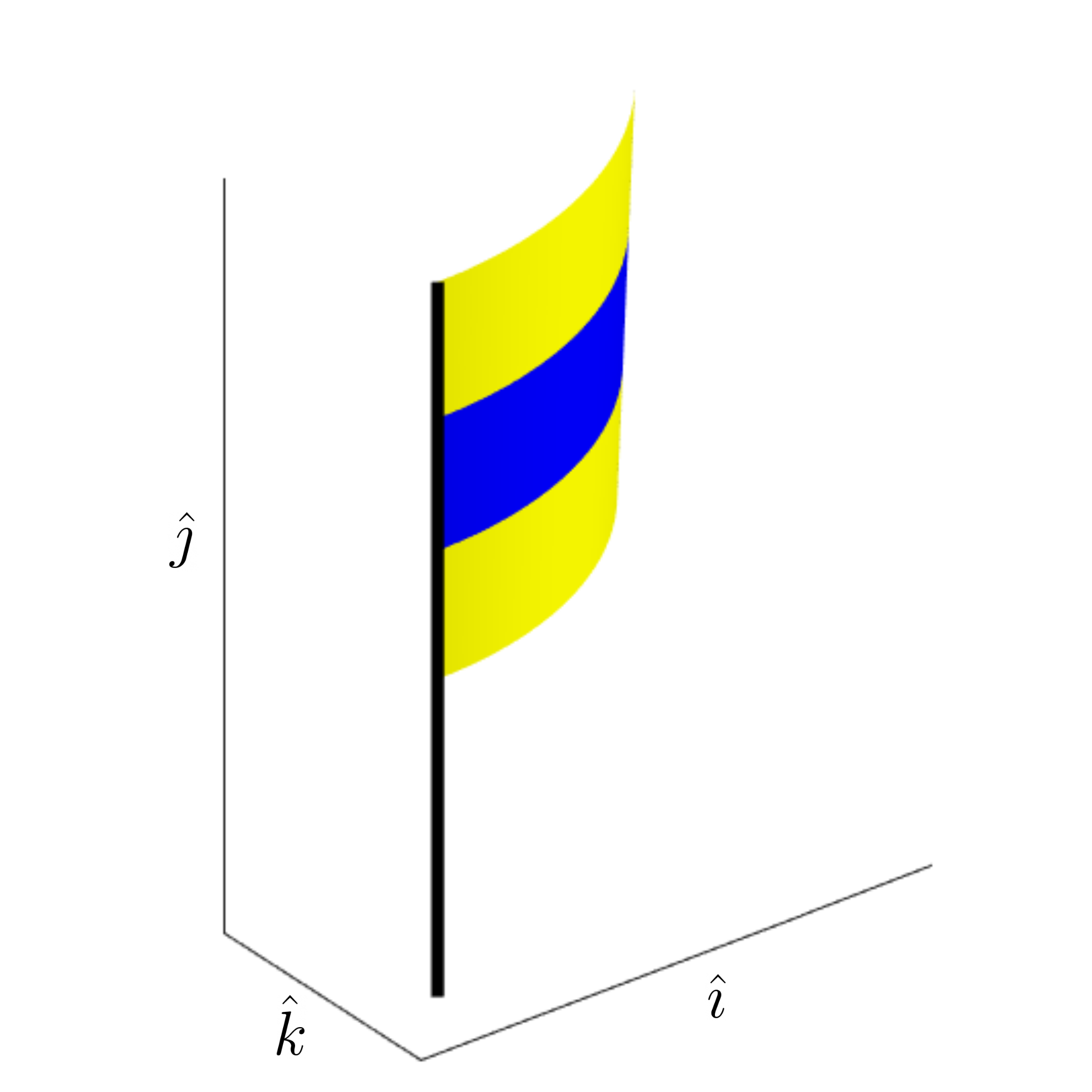}
				\includegraphics[width=.3\textwidth]{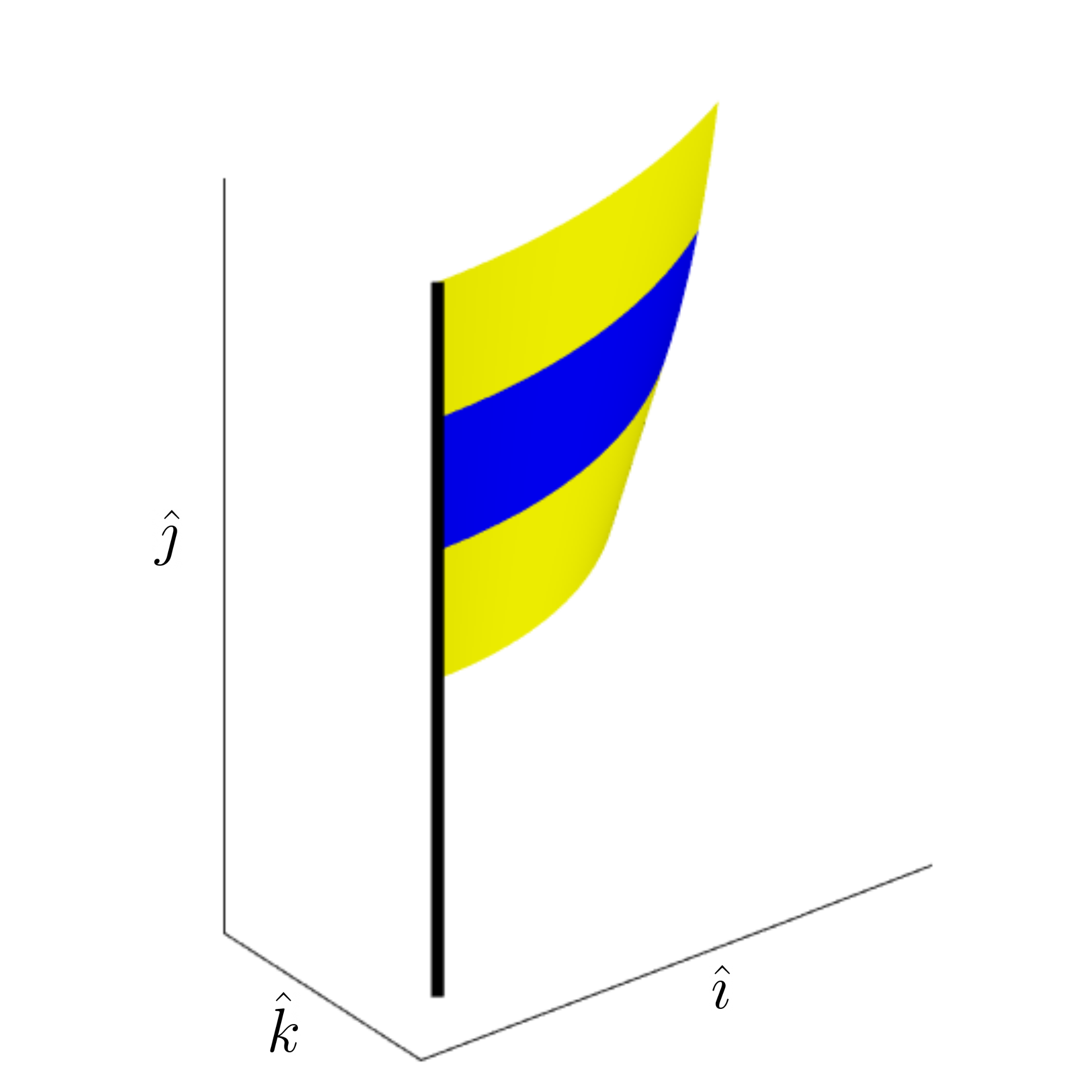}
						\includegraphics[width=.3\textwidth]{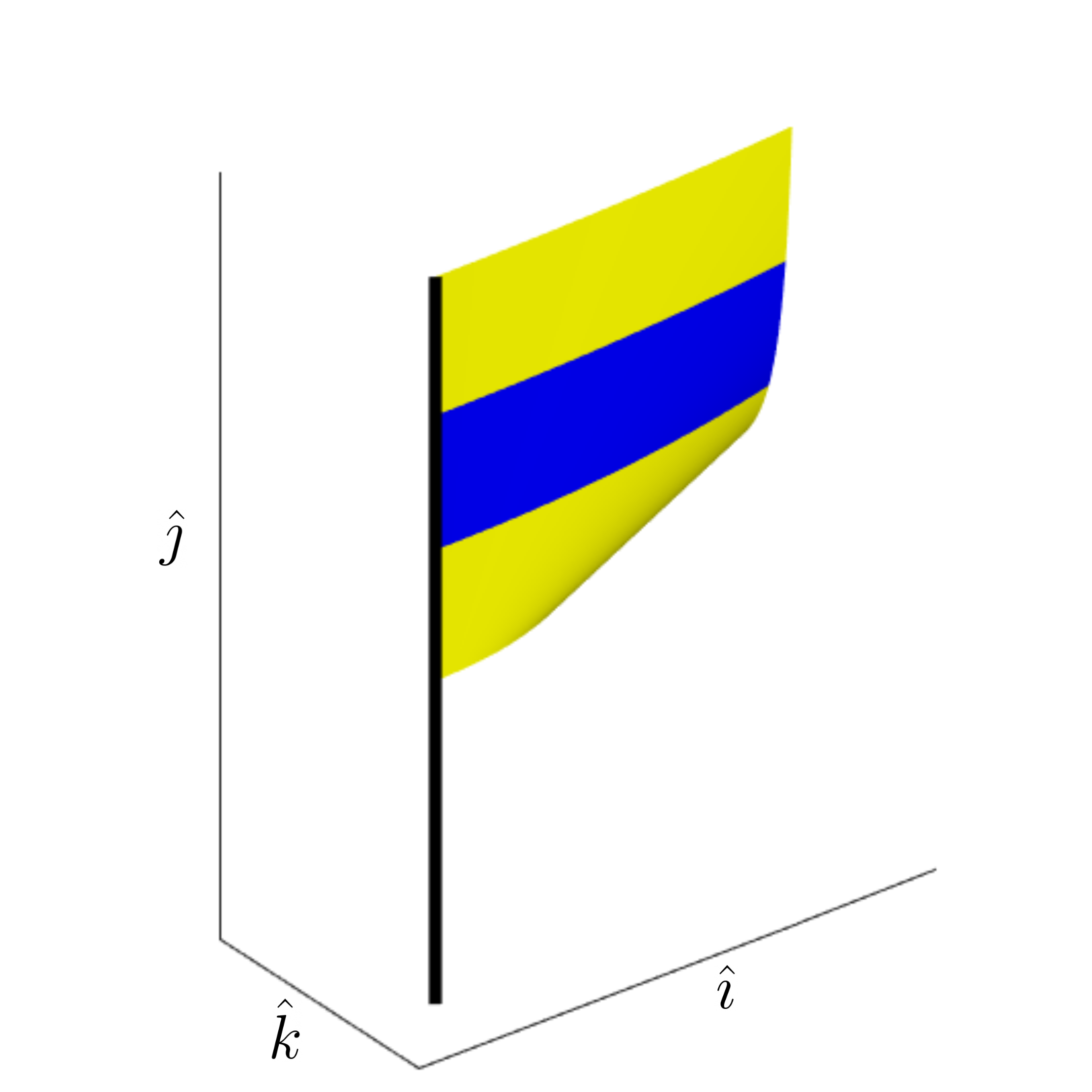}
\end{center}
\caption{Three examples of isometric immersions from the plane in $\mathbb R^3$ (flags). The flag pole is added for illustrative purposes only. The immersions have been constructed using the explicit characterization of Theorem~\ref{flagsolutionthm}.}
\label{fig:flagexamples}
\end{figure}

Our model for a flag is a function $\mathbf{r}\colon [0,1]^2 \to \mathbb{R}^3$ satisfying the isometry conditions $\lvert \mathbf{r}_u\rvert = \lvert \mathbf{r}_v\rvert = 1$ and $\mathbf{r}_u\cdot \mathbf{r}_v = 0$, along with the flagpole conditions $\mathbf{r}(0,v) = (0,v,0){:=v\hat{\mathbf{\jmath}}}$ and the horizontality condition $\mathbf{r}_u(0,v) ={(1,0,0):=} \hat{\mathbf{\imath}}$; in other words, the flag is fastened to the pole at $u=0$ with fasteners along the $\hat{\mathbf{\imath}}$ direction.
 {Here $\hat{\mathbf{\imath}}$, $\hat{\mathbf{\jmath}}$ and $\hat{k}$ denote a unit basis of $\mathbb R^3$.}

The first task is to classify maps satisfying these conditions. It is well-known that locally such immersions are determined by two functions of one variable; see e.g., the classic textbook by do Carmo~\cite{do2016differential}. {However, we need a global representation of these immersions to study the dynamics.} In our situation it turns out that either the curve along the bottom edge of the flag or the curve along the top edge of the flag has a special role, and we will distinguish these two cases by referring to them as upturned (downturned, resp.) flags. In the case that both of these curves admit this special role, we call it a balanced flag; see Section~\ref{sec:flagspace} for further details. Since the analysis for upturned and downturned flags is entirely equivalent, we will only focus on the case of upturned flags and only comment on the minor differences in Section~\ref{sec:downturned}. Our first main result, Theorem~\ref{flagsolutionthm}, gives a full characterization of all upturned flags in terms of two functions of one variable. The main difference between this and the classical results is that we get a global characterization on the whole square, which requires an analysis of crossing characteristics.

In essence, an upturned flag is determined by the space curve that traces out its bottom edge, which is in turn determined by its torsion and curvature functions, as long as these functions satisfy some constraint to keep the asymptotic lines from crossing. As the set of these functions is an open subset of a Banach space, this characterization provides us at the same time with a manifold structure for the space of upturned flags; see Theorem~\ref{flag_manifold}. It seems natural to consider the space of regular, upturned flags as a submanifold of the space of all $C^2$-surfaces. Since there happens to be a loss of derivatives, similarly to Nash's original investigations of the space of isometric immersions, this result seems unfortunately not true. This forces us to work in the smooth category, where we obtain a submanifold result using the inverse function theorem of Nash-Moser~\cite{hamilton1982inverse}, cf. Theorem~\ref{thm:submanifold}.

In the second part we study the natural kinetic energy metric on the space of flags, which allows us to model the motion of a flag as a geodesic curve with respect to this Riemannian structure. Towards this aim we then calculate the geodesic equation, which is obtained from the general principle that geodesics in a submanifold of a flat space satisfy the condition that the acceleration is normal to the submanifold. Deriving these equations  turns out to be the most involved computation of this article.
The complications illustrate the difficulty in modeling cloth or other unstretchable materials: isometric immersions are relatively rigid, but have some flexibility in special cases. This flexibility depends very much on the precise {boundary} conditions, however. In the final section, we present some preliminary numerical experiments for the geodesic boundary value problem that use the expression of the kinetic energy in terms of the two generating functions as derived in Section~\ref{sec:kinetic_energy}.

\subsection{Future directions}
In future work we plan to continue this line of research in several directions.
First, it would be of particular interest to obtain similar results for a more general class of isometric immersions.
The difficulty with actual isometrically embedded surfaces in $\mathbb{R}^3$ is that those without boundary must be fairly rigid, while those with boundary generate very complicated boundary conditions. As a first step to understanding these spaces it might be worth considering surfaces in $\mathbb{R}^4$, as it is much easier to isometrically embed them in this higher dimensional space; for example the tangent space at the standard Clifford torus can be written in terms of functions of two variables, not the single-variable functions that this quasi-rigidity gives us. Hence the theory will be more similar to that for the motion of inextensible closed curves in $\mathbb{R}^2$ or $\mathbb{R}^3$. Although the practical applications are obviously fewer, it would be an interesting space to study and perhaps reveal some information about the geometry of isometric immersions. Second, from an application-oriented point of view, we would like to use our geometric framework for the actual modeling of fabric; see for example \cite{ferreira2009shape}, \cite{chhatkuli2014non}, and \cite{WARDETZKY2007499} and references therein for discussions of current numerical methods for modeling fabric. This will require us to develop a comprehensive numerical framework for the calculation of both the geodesic initial and boundary value problem.

{
Furthermore, to model the movement of a real flag, one would want to incorporate the external force of gravity and the effect of wind (say, a uniform breeze in a fixed direction for simplicity).
Indeed there is a rich literature on modeling the interaction of flags with the surrounding fluid (wind), see eg.~\cite{zhang2000flexible}, \cite{1968JPSJ...24..392T}, \cite{shelley2011flapping}, \cite{argentina2005fluid} or \cite{fitt2001unsteady}. In our setup we would like to view this interaction with the external forces as an additional (potential) energy term, that we could then add to the total Lagrangian. 
While it is clear how to describe the potential energy due to gravity, it is far less straightforward how to incorporate even a simple model of wind. In the present article we describe these considerations briefly in Section~\ref{sec:gravity}, but a detailed study including the derivation of the resulting evolution equations is left open for future work.}

 Finally we would want a local existence theory for solutions of the geodesic equation in the space of flags. This already has major complications in the simpler case of inextensible curves (to which the flag equations reduce when nothing depends on the second spatial variable), as in the third author's paper~\cite{preston2011motion}. There we had a single wave equation with a tension determined by the solution of a second-order spatial ODE boundary value problem for each fixed time; here we have two coupled wave equations and a system of six first-order ODEs for each fixed time to determine the tensions. We have not attempted to address the local existence theory, since merely writing down the equations presents enough difficulty.
 
 {
 \subsection{Acknowledgements}
 The authors are grateful to Cy Maor for helpful discussions regarding the regularity assumptions on the space of isometric immersions. In addition the authors want to express their gratitude to the anonymous referee for their interesting comments, which greatly helped to improve the present article.  In particular, the suggestions on how to include the effects of wind into the present model, cf. Section~\ref{sec:gravity}, opened up an interesting avenue for future research. 
}
 
 {MB was partially supported
by NSF-grants 1912037 and 1953244 and by FWF-grant FWF-P 35813-N. SCP was partially supported by a Simons Foundation Collaboration Grant for Mathematicians no. 318969 and by a PSC-CUNY Award, jointly funded by The Professional Staff Congress and The City University of New York.}

 {Data sharing is not applicable to this article as no datasets were generated or analysed during the current study.}

\section{The space of flags}\label{sec:flagspace}
In this section we will introduce the basic notation of a flag (Definitions~\ref{flag}, \ref{regupflagdef} and \ref{regdownflagdef}) and show that any flag can be characterized uniquely  by two functions of one variable (Theorem~\ref{flagsolutionthm}). This will allow us to equip the space of all flags with a manifold structure (Theorem~\ref{flag_manifold}) and characterize its tangent space (Proposition~\ref{tangentspaceprop}). In addition we will show that we can view the space of smooth, regular flags as a submanifold of all smooth surfaces, cf. Theorem~\ref{thm:submanifold}. This will require us to use the Nash-Moser implicit function theorem.

We start by introducing the basic definition. \begin{definition}\label{flag}
{Let $\flag\in C^2([0,1]\times[0,1],\mathbb{R}^3)$.} We call $\flag$ a \emph{flag} if it is an isometric embedding of the square into $\mathbb{R}^3$
such that $\flag(0,v) = v\hat{\mathbf{\jmath}} $ and $\flag_u(0,v) = \hat{\mathbf{\imath}}$ for $v\in [0,1]$, {where $\hat{\mathbf{\imath}}, \hat{\mathbf{\jmath}}, \hat{k}$ is the standard basis of $\mathbb{R}^3$.}
We then have
\begin{equation}\label{orthonormal}
\flag_u\cdot \flag_u=1, \qquad \flag_u\cdot \flag_v = 0, \qquad \flag_v\cdot \flag_v = 1.
\end{equation}
\end{definition}

{
\begin{remark}[Nonlinear bending theory vs constrained membrane theory]
Note that one can make sense of the notion of a flag (isometric immersion, resp.) if the representing function $\flag$ is only in $C^1$. In the above definition we nevertheless require $\flag$ to be of class $C^2$. This regularity assumption has significant effects on the corresponding model of a flag:   by the celebrated results of Nash and Kuiper the space of $C^1$ isometric embeddings is dense in the space of short maps~\cite{nash1954c1,chern1952some}, whereas a similar statement is clearly not true for the space of $C^2$ isometric embeddings, cf.~\cite{conti2012h}. From a perspective of nonlinear elasticity theory, the $C^2$ assumption puts us in the realm of 
nonlinear bending theory, whereas the $C^1$ assumption can be viewed in the context of constrained membrane theory, cf. the seminal paper by Friesecke, James and M\"uller~\cite[Theorem 1, case (ii) and (iii)]{friesecke2006hierarchy}. To be more precise, the  appropriate modeling space for nonlinear bending theory is the Sobolev space of all $W^{2,2}$-isometric immersions; by a result of Pakzad~\cite{Pak04} the space of smooth ($C^2$, resp.) is dense in the space of $W^{2,2}$ immersions, which connects it to the regularity assumption of the present article.  From a mathematical point of view the $C^2$-assumption allows us to work with a continuous second fundamental, which serves as the basis for our chart construction. With exception of the submanifold result, we believe that all the constructions and results could be directly generalized to the $W^{2,2}$ category but beyond that, we suspect that entirely different techniques would be necessary.
From a modeling point of view the space of $C^1$ immersions (corresponding to constrained membrane theory) would allow for kink-like singularities that might appear while folding or crumbling a sheet of paper. We believe that excluding such irregularities is a reasonable assumption in the context of a cloth-made flag, which justifies the regularity requirements of the present article.  
\end{remark}
}

{
\begin{remark}[Rectangular flags]
In this article we  restrict ourself to flags that have a square shape. In practice, most flags are rather of a rectangular shape. By adding  additional parameters to the definition of the space of flags, one could easily extend the analysis of this article to this situation,  i.e., by considering isometric immersions of the form $\flag\colon[0,U]\times[0,V]\to\mathbb{R}^3$ with $U,V>0$. As the notation in the present article is already somewhat cumbersome, we have refrained from introducing these additional parameters. We want to emphasize that all the results of the article are also true in that situation; this has also been used in the illustrative example in Figure~\ref{fig:flagexamples}.
\end{remark}}

Letting   $\mathbf{N}=\flag_u\times \flag_v$ denote the normal vector, the fields $\{\flag_u, \flag_v, \mathbf{N}\}$ form a
convenient orthonormal basis. We compute that
\begin{equation}\label{efg}
\flag_{uu} = e(u,v) \mathbf{N}, \qquad \flag_{uv} = f(u,v) \mathbf{N}, \qquad \flag_{vv} = g(u,v) \mathbf{N},
\end{equation}
where the functions $e,f,g$ satisfy the Gauss-Codazzi and Codazzi-Mainardi equations
\begin{equation}\label{gausscodazzi}
eg - f^2 = 0, \qquad e_v = f_u, \qquad f_v = g_u.
\end{equation}
The second fundamental form $\big(\begin{smallmatrix} e & f \\ f & g\end{smallmatrix}\big)$ always has $0$ as an eigenvalue with eigenvector
$\big(\begin{smallmatrix} -f \\ e \end{smallmatrix}\big)$; the other eigenvalue is the mean curvature $e+g$. See for example \cite{Spi1979} or \cite{do2016differential}.

\begin{figure}[!ht]
\begin{center}
\includegraphics[scale=0.35]{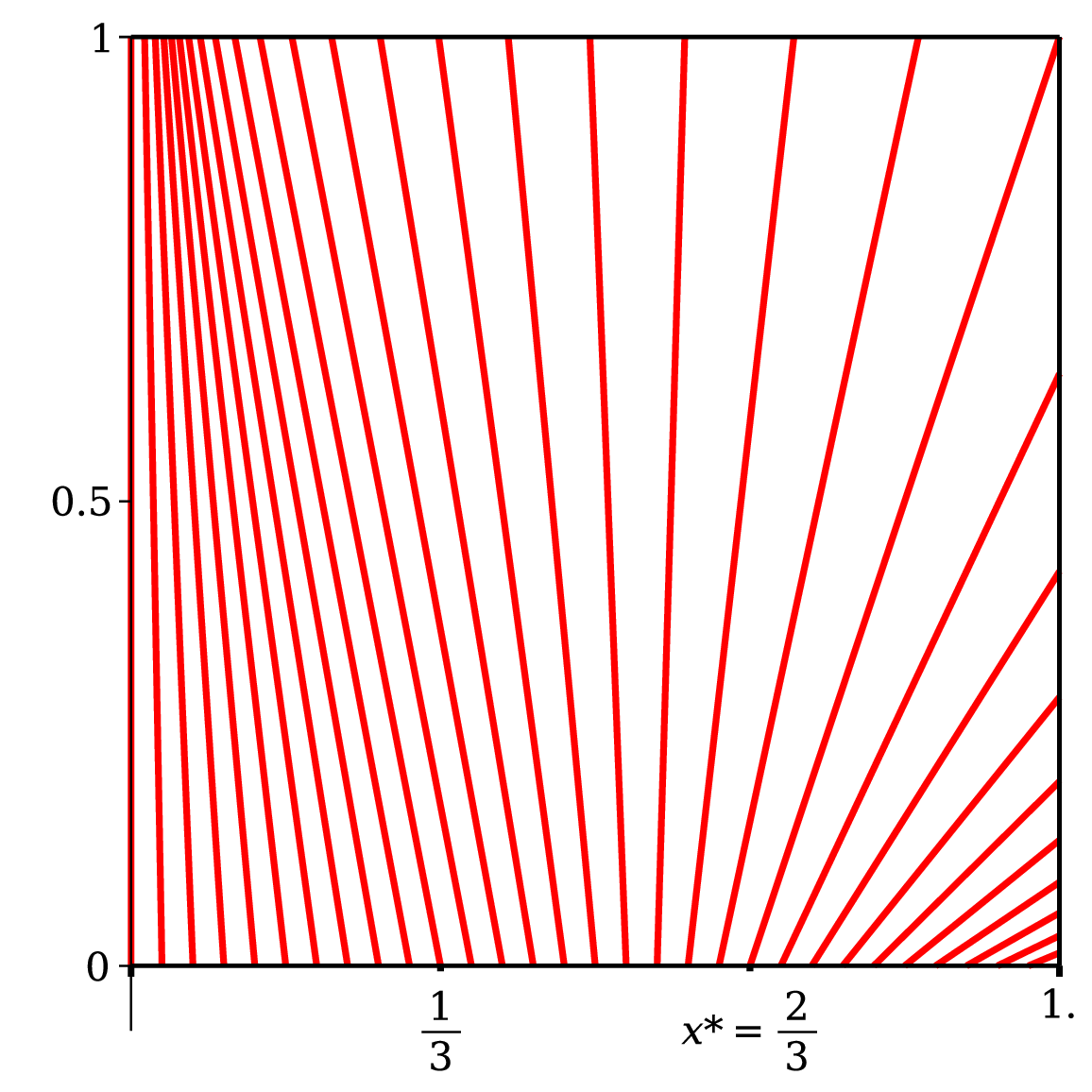} \qquad
\includegraphics[scale=0.35]{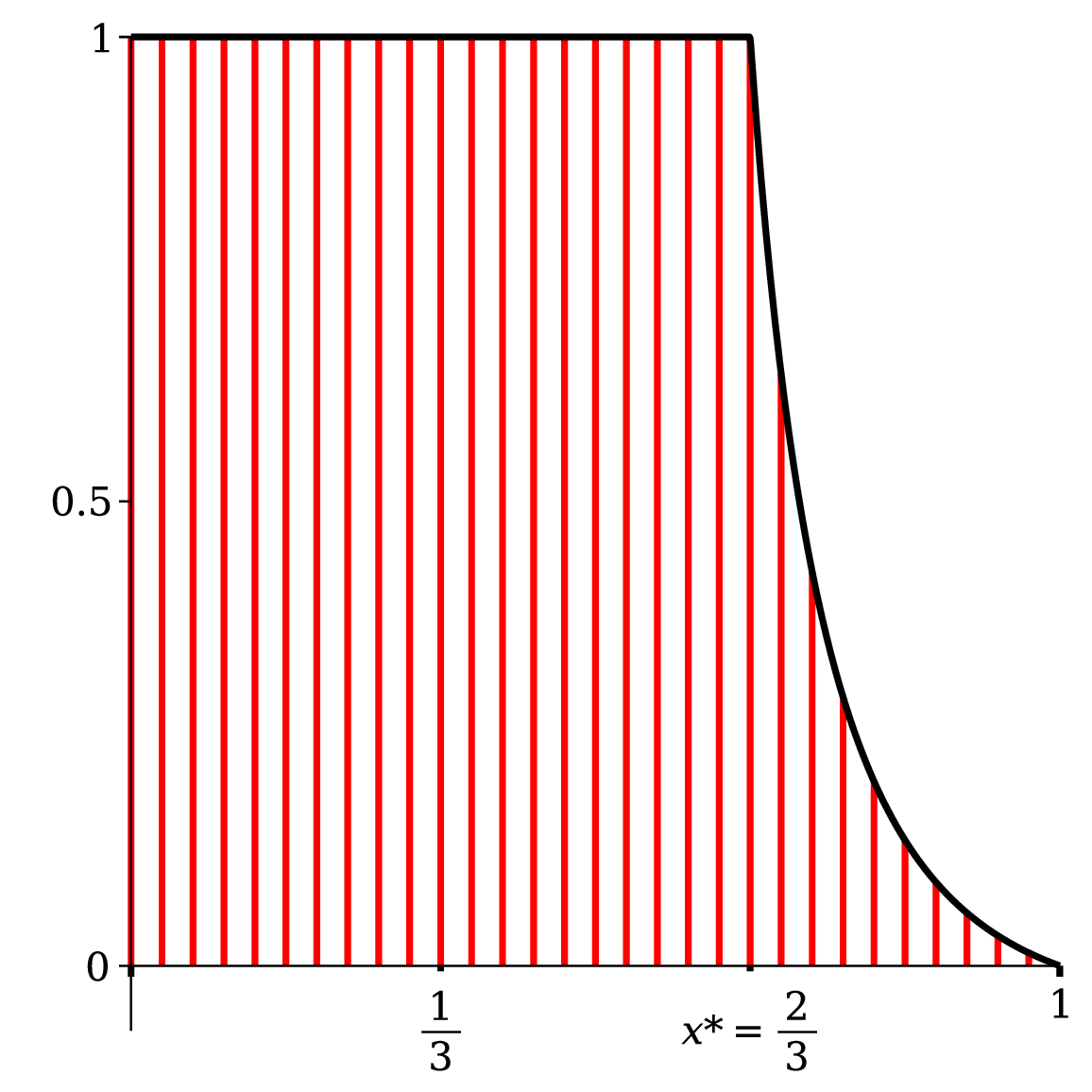}
\end{center}
\caption{On the left, we plot asymptotic lines on the unit square in the $(u,v)$ space for the function $\alpha_b(x) = -\tfrac{1}{2}x - 2x^2 + \tfrac{21}{4} x^3$, for which $x^*=\frac{2}{3}$. On the right, we plot these same lines in the new coordinates $(x,y)$ given by \eqref{bottomheavycoords}, for $y\le \gamma_b(x)$ as in \eqref{gammabdef}. The condition in Theorem \ref{gammaclassificationbottom} ensures that the right side of this shape is always a function that strictly decreases from $(x^*,1)$ to $(1,0)$.}
\label{fig:northdakota}
\end{figure}

We assume for simplicity that there is no open set of flat points where $e=f=g=0$, although we allow such points to occur. By a theorem, originally due to Pogorelov and Hartman-Wintner, as quoted in Ushakov~\cite{ushakov1996parameterisation}, there exists, for each point of the square a unique line segment through the point extending to the boundary of the flag.\footnote{We are unable to find the original statement and proof of this theorem, so we will assume it as a condition for nondegenerate flags.} The flagpole is one of these lines, and thus no other line can be horizontal (or it would cross the flagpole).
Consequently, we will consider the following three mutually exclusive cases:
\begin{itemize}
\item {\bf Upturned flag:} the {asymptotic} line through the upper right corner $(1,1)$ passes through the bottom edge at a point $(x^*,0)$ for $0< x^*<1$.
\item {\bf Downturned flag:} the {asymptotic} line through the lower right corner $(1,0)$ passes through the top edge at $(x^*,1)$ for $0< x^*<1$.
\item {\bf Balanced flag:} a single asymptotic line passes through both right corners.
\end{itemize}
In the upturned case, asymptotic lines starting on the bottom edge will hit every point of both the top edge and the right side, so we parameterize everything along this bottom edge; the asymptotic line through the bottom right corner reduces to a point. Similarly in the downturned case, asymptotic lines through points of the top edge will pass through both the bottom edge and the right side, and we parameterize along the top edge. In the balanced case, we could use either parameterization along the top edge or bottom edge. We refer to Figure~\ref{fig:northdakota} for an example of an upturned flag, which {illustrates} the above construction.

We start by defining a coordinate transformation, which will be central in the remainder of the article. We use $x$ as the parameter along the bottom or top edge (depending on whether the flag is upturned or downturned), and $y$ as the parameter along the lines.
\begin{itemize}
\item For an upturned flag, let $\alpha_b(x)$ be the reciprocal of the slope of the asymptotic line through $(x,0)$, so that the asymptotic line is given by
    \begin{equation}\label{bottomheavycoords}
    u=x+\alpha_b(x)y, \qquad v=y
     \end{equation}
     for $0\le x\le x^*$, $0\le y\le 1$ (to fill out the top edge) and $x^*\le x\le 1$, $0\le y\le \frac{1-x}{\alpha_b(x)}$ (to fill out the right edge). {We denote the graph of $y$ as a function of $x$ by \begin{equation}\label{gammabdef}\gamma_b(x)=
     \begin{cases} 1 & 0\le x\le x^* \\
     \frac{1-x}{\alpha_b(x)} & x^*< x\le 1.\end{cases}\end{equation}}
\item For a downturned flag, let $\alpha_t(x)$ be the reciprocal of the slope of the asymptotic line through $(x,1)$, so the line is given by \begin{equation}\label{topheavycoords}
    u=x-\alpha_t(x)y, \qquad v=1-y
     \end{equation}
     for $0\le x\le x^*$ and $0\le y\le 1$, or $x^*\le x\le 1$ and $0\le y\le -\frac{1-x}{\alpha_t(x)}$. {Similarly we have a piecewise-defined function $\gamma_t(x)$ which traces the bottom edge.}
\item For a balanced flag we may use either representation, where $x^*=1$.
\end{itemize}

Note that when the asymptotic line passes through both the top and bottom of the flag (rather than the right side), the line given through the bottom \eqref{bottomheavycoords} for some $x$ passes through the point $(z,1)$ for some $z$, and we have $z=x+\alpha_b(x)$. The slope of this line is $\alpha_b(x)$, but it is also equal to $\alpha_t(z)$, and therefore we must have
\begin{equation}\label{alphabalphat}
\alpha_b(x) = \alpha_t\big(x+\alpha_b(x)\big), \text{ if $\alpha_b(x)+x<1$.}
\end{equation}
Similarly we have
\begin{equation}\label{alphatalphab}
\alpha_t(x)=\alpha_b\big(x-\alpha_t(x)\big), \text{ if $x-\alpha_t(x)<1$.}
\end{equation}
In the upturned case there is an $x^*$ with $\alpha_b(x^*)+x^*=1$, and here we will have $\alpha_t(1) = \alpha_b(x^*)=1-x^*$; in particular $\alpha_t(1)>0$. Similarly in the downturned case there is an $x_*$ with $x_*-\alpha_t(x_*)=1$, and in this case we will have $\alpha_b(1) = \alpha_t(x_*)<0$. In the balanced case we have $\alpha_b(1)=\alpha_t(1)=0$. We can thus distinguish the three cases in terms of the single number $\alpha_b(1)$: if positive the flag is upturned, if negative the flag is downturned, and if zero the flag is balanced. Similarly we could do the same with $\alpha_t$.

At points where $e\ne 0$ we may define the ratio $\phi(u,v) = f(u,v)/e(u,v)$. Since $\big(\begin{smallmatrix}  -f \\ e\end{smallmatrix}\big)$ is the vector in the direction of the nullspace, the functions $\alpha_b$ and $\alpha_t$ are related to $\phi$ by
$$ \phi(x,0) = -\alpha_b(x), \qquad \phi(x,1) = -\alpha_t(x).$$
The following lemma explains the usefulness of our new coordinates $(x,y)$.

\begin{lemma}\label{burgers}
{Let $\flag\in C^2([0,1]\times[0,1],\mathbb R^3)$ and let $e$ and $f$ be given by~\eqref{efg}.}
At points where $e\ne 0$, the function $\phi = f/e$ satisfies the inviscid Burgers' equation $\phi_v = \phi \phi_u$, and thus in $(x,y)$ coordinates we have $\phi_y=0$.
\end{lemma}

\begin{proof}
By the Gauss-Codazzi equation \eqref{gausscodazzi}, the zero-curvature condition $eg=f^2$ becomes $f=\phi e$ and $g=\phi^2e$. The Codazzi-Mainardi equations then imply that
$$ \frac{\partial}{\partial u}(\phi^2 e) - \phi \,\frac{\partial}{\partial v}(\phi e) = \phi \left( \frac{\partial}{\partial u}(\phi e) - \phi \,\frac{\partial}{\partial v}(e)\right),$$
which reduces to $e \phi\phi_u = e\phi_v$, yielding the inviscid Burgers' equation, which is easily solved by the method of characteristics (the characteristics being precisely the asymptotic lines).

We conclude that if the initial data is given along the bottom edge (as in the upturned case), then
$\phi$ satisfies the implicit equation
\begin{equation}\label{bottomsoln}
\phi(u,v) = -\alpha_b\big( u - v \phi(u,v)\big),
\end{equation}
while if the data is given along the top edge, then $\phi$ satisfies
\begin{equation}\label{topsoln}
\phi(u,v) = -\alpha_t\big(u+(1-v)\phi(u,v)\big).\qedhere
\end{equation}
 \end{proof}

The next corollary simply follows from the coordinate transformations \eqref{bottomheavycoords}--\eqref{topheavycoords} and the solution formulas \eqref{bottomsoln}--\eqref{topsoln}.

\begin{corollary}\label{phisoln}
In the $(x,y)$ coordinates given in the upturned case by \eqref{bottomheavycoords}, we may write $\phi(u,v) = -\alpha_b(x)$. Similarly in the downturned  case with coordinates~\eqref{topheavycoords}, we may write $\phi(u,v) = -\alpha_t(x)$. In he balanced case we can use either.
\end{corollary}
From here on, we will focus solely on the upturned case; all the following results also hold in the two other cases and we will comment on the small differences in the end of the section. In what follows we will skip the subscript $b$ in functions such as $\alpha$.

\subsection{A characterization for upturned flags}
In this part we will show in the upturned case that the function $\alpha$ determining the asymptotic lines, together with a curvature function $\kappa(x):= e(x,0)$ specified along the bottom edge, completely determine the flag. In fact this is essentially a statement that the upturned flag is completely determined by the unit-speed curve $\boldsymbol{\eta}(x) = \flag(x,0)$ along the bottom of it.

We start by characterizing all functions $\alpha$ that create an upturned flag.  The central ingredient for this result is the observation that the asymptotic curves are characteristics of a homogeneous quasilinear PDE and that the solutions of this PDE can be differentiable iff the characteristics do not cross. \begin{theorem}\label{gammaclassificationbottom}
{A function $\alpha\in  C^1([0,1],\mathbb{R})$} with $\alpha(0)=0$ and $\alpha(1)>0$ generates a family of nonintersecting asymptotic lines in the upturned flag case, with the coordinate transformation \eqref{bottomheavycoords} forming a global diffeomorphism on the square, iff it satisfies the condition
 \begin{equation}\label{lambdapositive}
 \lambda(x)>0  \text{ for all $x\in[0,1]$,}
 \end{equation}
 where
\begin{equation}\label{lambdagammadef}
\lambda(x) = 1+\alpha'(x)\gamma(x), \qquad
\gamma(x) = \begin{cases} 1 & \alpha(x)\le 1-x, \\
\frac{1-x}{\alpha(x)} & \alpha(x)>1-x.
\end{cases}
\end{equation}
\end{theorem}

\begin{proof}
First suppose $\alpha\colon[0,1]\to\mathbb{R}$ generates a family of nonintersecting lines filling up the square. The lines are given in the square $(u,v)\in [0,1]^2$ by the parameterization
$u = x+\alpha(x)v$, for $0\le x\le 1$. For the upturned flag, there is some $x^*$ such that the line passes through the top right corner $(1,1)$, so that we have $1-x^*=\alpha(x^*)$. We will refer to this special line as the corner-bending line. Again see Figure \ref{fig:northdakota} for the visual interpretation.

For all lines with $x\le x^*$, to the left of the corner-bending line, the vertical parameter $v$ goes from $0$ to $1$ as $u$ goes from $x$ to $x+\alpha(x)$. In particular, we know this line does not reach the right corner, so that $x+\alpha(x)<1$ for $x<x^*$. On the other hand, to the right of the corner-bending line, the asymptotic line crosses the right side of the square $u=1$ before $v$ reaches $1$; in fact at height $v=\frac{1-x}{\alpha(x)}>0$.
In particular we see that $\alpha(x)>1-x$ for $x>x^*$. Thus $\gamma$ defined by \eqref{lambdagammadef}, the largest value of $v$ such that the line segment remains in the square, satisfies $$\gamma(x) = \begin{cases} 1 & x\le x^*, \\
\frac{1-x}{\alpha(x)} & x>x^*.
\end{cases}$$
Hence we have a parameterization of the unit square $(u,v)\in [0,1]^2$ by the map
\begin{equation}\label{xytouv}
(x,y)\to (u,v) = (x+\alpha(x)y, y), \qquad 0\le x\le 1, \qquad 0\le y\le \gamma(x).
\end{equation}
This parameterization is an invertible diffeomorphism on the entire square if and only if its Jacobian determinant is everywhere positive. That condition on the Jacobian is clearly
\begin{equation}\label{jacobian}
\text{Jac} = 1+\alpha'(x)y>0 \qquad \text{for $0\le y\le \gamma(x)$},
\end{equation}
which is equivalent to $\lambda(x)>0$ by the definition \eqref{lambdagammadef}.

Conversely suppose a function $\alpha$ which generates
$\gamma$ and $\lambda$ by formula \eqref{lambdagammadef} has $\lambda(x)>0$ for all $x\in [0,1]$. We want to show the parameterization \eqref{xytouv} fills up the square $(u,v)\in [0,1]^2$. For each fixed $(u,v)\in [0,1]^2$, define $F_{(u,v)}\colon [0,1]\to \mathbb{R}$ by
\begin{equation}\label{Fdef}
F_{(u,v)}(x) = x + v\alpha(x) - u.
\end{equation}
We want to show that for each $(u,v)$ there is a unique $x$ such that $F_{(u,v)}(x)=0$, which will imply $(u,v)$ is reached by the parameterized line. Clearly if $u=0$, then $x=0$ is such a point, since $F_{(u,v)}(0) = 0+v\alpha(0)=0$; while if $v=0$, then $x=u$ is obviously the unique such point. So we will assume $u>0$ and $v>0$ in what follows.

Since $\alpha(0)=0$ and $\alpha(1)>0$, we have
$$ F_{(u,v)}(0) = -u<0 \qquad \text{and}\qquad F_{(u,v)}(1)=v\alpha(1)>0,$$
so there is always at least one such point. We now want to show uniqueness, which we will do by showing that $F_{(u,v)}'(x)>0$ whenever $F_{(u,v)}(x)=0$.
So suppose $F_{(u,v)}(x)=0$, and compare $\alpha(x)$ to $1-x$.
\begin{itemize}
\item If $\alpha(x)\le 1-x$, then $\gamma(x)=1$ and $\lambda(x)=1+\alpha'(x)$ by \eqref{lambdagammadef}. The inequality \eqref{lambdapositive} implies that $\alpha'(x)>-1$, so that
    $F_{(u,v)}'(x) = 1+v\alpha'(x) > 1-v\ge 0$.
\item If $\alpha(x)>1-x$, then $\gamma(x)=\frac{1-x}{\alpha(x)}$, and
$\lambda(x)=1+\frac{\alpha'(x)(1-x)}{\alpha(x)}$ by \eqref{lambdagammadef}.
In addition, since $F_{(u,v)}(x)=0$, we get
$$ v\alpha(x) = u-x,$$
so that
\begin{align*}
F_{(u,v)}'(x) &= 1 + v\alpha'(x) = 1 + \frac{v\alpha(x) (\lambda(x)-1)}{1-x} \\
&> 1 - \frac{v\alpha(x)}{1-x} = 1 - \frac{u-x}{1-x} = \frac{1-u}{1-x} \ge 0.
\end{align*}
\end{itemize}
Either way, $F_{(u,v)}'(x)>0$ whenever $F_{(u,v)}(x)=0$, so there is exactly one $x$ such that $F_{(u,v)}(x)=0$ for each $u>0$ and $v>0$.

In particular there is a unique $x^*$ satisfying $\alpha(x^*)=1-x^*$, and for $x<x^*$ we must have $\alpha(x)<1-x$, while for $x>x^*$ we must have $\alpha(x)>1-x$. Hence the formula \eqref{lambdagammadef} defining $\gamma$ becomes
$$\gamma(x) =\begin{cases} 1 & x\le x^*, \\
\frac{1-x}{\alpha(x)} & x>x^*.
\end{cases}$$
The map \eqref{bottomheavycoords} is thus a bijection for $0\le y\le \gamma(x)$ and $0\le x\le 1$ onto the square. Since its Jacobian determinant is given by $J(x,y) = 1+y\alpha'(x)$, which for each fixed $x\in [0,1]$ ranges from $1$ to $\lambda(x)$, we see this Jacobian determinant is positive, so the map \eqref{bottomheavycoords} is a global diffeomorphism.
\end{proof}

Based on the discussion above, we make the following definition for an upturned flag affixed to a vertical flagpole (pointing in the $\hat{\jmath}$ direction) along a horizontal edge (in the $\hat{\imath}$ direction). 
We require that $\flag$ is a $C^2$ function, so that $\flag_{uu}$, $\flag_{uv}$, and $\flag_{vv}$ are all $C^0$, but in addition we require that $\alpha$ is a $C^1$ function (which is not automatic, since it is the ratio of functions which are only a priori continuous).

\begin{definition}\label{regupflagdef}
A \emph{regular upturned flag} is a $C^2$ isometric embedding $\flag\colon [0,1]^2 \to \mathbb{R}^3$ satisfying the conditions
\begin{equation}\label{flagboundaryconditions}
\flag(0,v)=v \hat{\jmath}, \quad \flag_u(0,v) = \hat{\imath}, 
\quad \text{for all $v\in [0,1]$,}
\end{equation}
and such that there is a $C^1$ function $\alpha\colon [0,1]\to \mathbb{R}$ satisfying
$\flag_{uv}(x,0) = -\alpha(x) \flag_{uu}(x,0)$ as well as the conditions of Theorem \ref{gammaclassificationbottom}. That is, $\alpha(0)=0$, $\alpha(1)>0$, and
\begin{equation}\label{nondegenerateflagcond}
\max\{1-x, \alpha(x)\} + (1-x)\alpha'(x) > 0 \quad \text{for all $x\in [0,1]$.}
\end{equation}
\end{definition}

Above we saw that given an isometric nonsingular immersion of the square, much of the geometry is characterized by the properties of a single function $\alpha\colon [0,1]\to \mathbb{R}$. It remains to show that this function, together with a curvature function $\kappa(x):= e(x,0)$ specified along the bottom edge, completely determines the flag. This is essentially a statement that the upturned flag is completely determined by the unit-speed curve $\boldsymbol{\eta}(x) = \flag(x,0)$ along the bottom of it, since that curve is uniquely determined by its curvature and torsion via the classification theorem for curves. Here the curvature is the function $\kappa(x)$, while the torsion is $\tau(x) = -\alpha(x)\kappa(x)$.

We will not quite take the usual Frenet-Serret approach, since we want to allow the curvature to change sign, which allows us to reproduce the two-dimensional case where the curvature is signed. The ordinary Frenet-Serret theory assumes that the curvature is never zero, so that the normal is always well-defined; that is an issue for a general space curve, but not in this situation since we already have an orthonormal frame $\{\flag_u, \flag_v, \mathbf{N}\}$. This approach also has the advantage that it explicitly reconstructs both the curve and the flag from the curvature and torsion, via a system of ordinary differential equations for the spherical coordinates.

\begin{lemma}\label{frenetserretlemma}
If $\alpha\colon [0,1]\to \mathbb{R}$ is $C^1$ and $\kappa\colon [0,1]\to\mathbb{R}$ is $C^0$, then there is a unique $C^2$ unit-speed curve $\boldsymbol{\eta}\colon [0,1]\to\mathbb{R}^3$ satisfying the Frenet-Serret equations:
\begin{equation}\label{frenetserreteqs}
\boldsymbol{\eta}''(x) = \kappa(x) \normal(x), \qquad \normal'(x) = -\kappa(x)\big(\boldsymbol{\eta}'(x) + \alpha(x)\binormal(x)\big), \qquad \binormal'(x) = \kappa(x)\alpha(x) \normal(x),
\end{equation}
with $\{\tangent = \boldsymbol{\eta}', \normal,\binormal\}$ forming a $C^1$ oriented orthonormal basis for each $x$, and such that  $\boldsymbol{\eta}(0)=0$, $\boldsymbol{\eta}'(0)=\hat{\imath}$, and $\normal(0)=\hat{k}$.
\end{lemma}

\begin{proof}
Since $\boldsymbol{\eta}'$ is to be a unit vector field, we define spherical coordinates by
\begin{equation}\label{etae1def}
\boldsymbol{\eta}'(x) = \big( \cos{\theta(x)} \cos{\phi(x)}, \cos{\theta(x)} \sin{\phi(x)}, \sin{\theta(x)}\big).
\end{equation}
The condition $\boldsymbol{\eta}'(0)=(1,0,0)$ means that $\phi(0)=\theta(0)=0$. We suppose $\theta\in [0,\pi)$ and $\phi\in [0,2\pi)$.

Set $e_1(x) = \boldsymbol{\eta}'(x)$, and complete to an orthonormal basis $\{e_1,e_2,e_3\}$ via the formulas
\begin{equation}\label{e2e3def}
\begin{split}
e_2(x) &= \big( \sin{\theta(x)}\cos{\phi(x)}, \sin{\theta(x)}\sin{\phi(x)}, -\cos{\theta(x)}\big) \\
e_3(x) &= \big( -\sin{\phi(x)}, \cos{\phi(x)}, 0\big).
\end{split}
\end{equation}
Setting
\begin{equation}\label{normalbinormaldef}
\begin{split}
\normal(x) &= -\cos{\beta(x)} e_2(x) + \sin{\beta(x)} e_3(x), \\
\binormal(x) &= -\sin{\beta(x)} e_2(x) - \cos{\beta(x)} e_3(x),
\end{split}
\end{equation}
for some function $\beta\colon [0,1]\to\mathbb{R}$, we see that $\boldsymbol{\eta}'(x)\times \normal(x) = \binormal(x)$ for all $x$. The condition $\normal(0)=\hat{k}$ together with $e_2(0)=-\hat{k}$ and $e_3(0) = \hat{\jmath}$ implies that we must have $\beta(0)=0$.

Now consider the system
\begin{alignat}{3}
\theta'(x) &= \kappa(x) \cos{\beta(x)}, &\qquad \theta(0)&=0; \label{thetaeq} \\
\phi'(x) &= \kappa(x) \sec{\theta(x)} \sin{\beta(x)}, & \phi(0)&= 0; \label{phieq} \\
\beta'(x) &= \kappa(x)\big( \alpha(x)+\sin{\beta(x)} \tan{\theta(x)}\big), & \beta(0)&= 0.\label{betaeq}
\end{alignat}
There is a unique solution $\{\theta,\phi,\beta\}$ for $x$ close to zero, and the solutions are $C^1$ functions.

We have
$$ \boldsymbol{\eta}''(x) = -\theta'(x) e_2(x) + \cos{\theta(x)} \phi'(x) e_3(x) = \kappa(x) \normal(x)$$
using equations \eqref{thetaeq} and \eqref{phieq}. In addition these equations imply
\begin{align*}
\binormal'(x) &=
-\cos{\beta(x)} \beta'(x) e_2(x) + \sin{\beta(x)} \beta'(x) e_3(x)  - \sin{\beta(x)}\big( \theta'(x) e_1(x) + \sin{\theta(x)} \phi'(x) e_3(x)\big)\\
&\qquad\qquad+ \phi'(x) \cos{\beta(x)}\big( \cos{\theta(x)} e_1(x) + \sin{\theta(x)} e_2(x)\big) \\
&= \big(\phi'(x) \cos{\beta(x)} \cos{\theta(x)} - \theta'(x) \sin{\beta(x)} \big) e_1(x) \\
&\qquad\qquad + \big(\beta'(x) - \phi'(x) \sin{\theta(x)}\big) (-\cos{\beta(x)} e_2(x) + \sin{\beta(x)} e_3(x)\big) \\
&= \kappa(x) \alpha(x) \normal(x).
\end{align*}

Orthonormality of the frame $\{\tangent, \normal, \binormal\}$ then implies the remainder of the Frenet-Serret equations, that
$$ \normal'(x) = -\kappa(x) \big( \tangent(x) + \alpha(x) \binormal(x)\big).$$
Furthermore since the frame $\{\tangent, \normal, \binormal\}$ remains orthonormal, the components remain bounded, and the solution of the ODE system exists for all $x\in [0,1]$, not just locally.

If $\alpha$ and $\kappa$ are at least $C^0$ functions, then the solution of the system \eqref{thetaeq}--\eqref{betaeq} must be $C^1$. This implies in particular that $\boldsymbol{\eta}'$ is $C^1$, so that $\boldsymbol{\eta}$ is $C^2$.
\end{proof}

In the following theorem we demonstrate that all regular upturned flags are completely characterized by the continuous curvature function $\kappa$ and the continuously differentiable function $\alpha$, both specified along the bottom edge. The formula \eqref{regupflagformula} is called the \emph{rectifying developable} or the \emph{envelope of tangent planes} of the bottom edge $\boldsymbol{\eta}$. It is a well-known classical result (see e.g., Struik~\cite{struik1961lectures} or do Carmo~\cite{do2016differential}) that any such surface is a developable surface, i.e., that it is locally the image of an isometric immersion of the plane. The reason for the somewhat more complicated presentation here is that we obtain a \emph{global} description of the flag on the entire square, not merely a local representation of it. This allows us to also consider the set of all regular upturned flags as a topological space and a manifold, and study geodesic motion in it, as we shall do later. The parameterization we give here ends up being close to that of Izumiya et al.~\cite{izumiya1999rectifying}, who coined the term ``modified Darboux vector'' for the vector $\darboux(x)$ below.

\begin{theorem}\label{flagsolutionthm}
For any regular upturned flag $\flag\colon [0,1]^2\to\mathbb{R}^3$ as in Definition \ref{regupflagdef}, the bottom edge defined by $\boldsymbol{\eta}(x) = \flag(x,0)$ is a unit-speed curve, with (signed) curvature
\begin{equation}\label{bottomcurvature}
\kappa(x) = e(x,0) = \langle \flag_{uu}(x,0), \flag_u(x,0)\times \flag_v(x,0)\rangle
\end{equation}
and torsion given by $\tau(x) = -\alpha(x) \kappa(x)$. The function $\alpha$ satisfies the conditions $\alpha(0)=0$, with $\alpha(1)>0$ and the inequality \eqref{nondegenerateflagcond}.
The Frenet-Serret frame $\{\tangent, \normal, \binormal\}$ along $\boldsymbol{\eta}$ is given by
$$ \tangent(x) = \flag_u(x,0), \qquad \normal(x) = \mathbf{N}(x,0), \qquad \binormal(x) = -\flag_v(x,0).$$

Conversely, given any $C^1$ function $\alpha\colon [0,1]\to \mathbb{R}$ satisfying the conditions of Theorem \ref{gammaclassificationbottom} or equivalently Definition \ref{regupflagdef}, and any $C^0$ function $\kappa\colon [0,1]\to \mathbb{R}$, there is a unique regular upturned flag $\flag$ given by
\begin{equation}\label{regupflagformula}
\flag(u,v) = \boldsymbol{\eta}(x) - v \darboux(x), \qquad \darboux(x) = -\alpha(x) \boldsymbol{\eta}'(x)+\binormal(x),
\end{equation}
where $x$ is defined for each $u,v\in [0,1]^2$ to be the unique solution in $[0,1]$ of
\begin{equation}\label{xcondition}
x + \alpha(x)v = u,
\end{equation}
and $\darboux(x)$ is the modified Darboux vector of $\boldsymbol{\eta}$.
\end{theorem}

\begin{proof}
Given the upturned flag $\flag$, the fact that $\flag$ is an isometry implies that $\flag_u(x,0) = \boldsymbol{\eta}'(x)$ is a unit vector for all $x\in [0,1]$. We have $\boldsymbol{\eta}''(x) = \flag_{uu}(x,0) = e(x,0) \mathbf{N}(x,0)$ by equations \eqref{efg}, and so if we define $\normal(x) = \mathbf{N}(x,0)$ to be the normal field along the curve, then $\boldsymbol{\eta}''(x) = \kappa(x) \normal(x)$ with $\kappa(x)=e(x,0)$, which is the first of the Frenet-Serret equations \eqref{frenetserreteqs}.
Since $\flag_u \times \mathbf{N} = -\flag_v$ by definition of $\mathbf{N}$, and $\tangent\times \normal =\binormal$ by construction of the Frenet-Serret basis, we must have $\binormal(x) = -\flag_v(x,0)$. Furthermore we compute
$$ \binormal'(x) = -\flag_{uv}(x,0) = -f(x,0) \mathbf{N}(x,0) = \alpha(x) e(x,0) \normal(x)$$
using \eqref{efg} and equation \eqref{bottomsoln} from Lemma \ref{burgers}, which is the third of the Frenet-Serret equations \eqref{frenetserreteqs}.
As in the proof of Lemma \ref{frenetserretlemma}, the formula for the derivative of $\normal$ is the second Frenet-Serret equation.

Now we consider the converse, supposing that $\alpha$ is $C^1$ and $\kappa$ is a given $C^0$ function satisfying $\alpha(0)=\kappa(0)=0$, $\alpha(1)>0$, and the inequality \eqref{nondegenerateflagcond}. Using Lemma \ref{frenetserretlemma}, we construct the unique curve $\boldsymbol{\eta}\colon [0,1]\to\mathbb{R}^3$ along with its orthonormal Frenet-Serret frame $\{\tangent=\boldsymbol{\eta}', \normal,\binormal\}$,
subject to the conditions $\boldsymbol{\eta}(0)=0$, $\boldsymbol{\eta}'(0) = \hat{\imath}$, and $\normal(0) = \hat{k}$. This curve $\boldsymbol{\eta}$ is $C^2$, and its Frenet-Serret frame is $C^1$. Then we define a surface by the parameterization \eqref{regupflagformula}. The fact that the function $(u,v)\mapsto x$ given by \eqref{xcondition} is well-defined and continuously differentiable is a consequence of our regularity definition and Theorem \ref{gammaclassificationbottom}.

We first show that this parameterized surface is a $C^2$ isometric immersion, and to do this we compute $\flag_u$ and $\flag_v$ and show that these are $C^1$ and orthonormal for all $(u,v)\in [0,1]^2$ as in \eqref{orthonormal}.

First we compute the derivatives of $x(u,v)$ implicitly from \eqref{xcondition}, which gives
\begin{equation}\label{xpartials}
\frac{\partial x}{\partial u}  = \frac{1}{1+v \alpha'(x)}, \qquad \frac{\partial x}{\partial v} = -\frac{ \alpha(x)}{1+v\alpha'(x)},
\end{equation}
and the fact that the denominators are always positive and well-defined for all $v\in [0,\gamma(x)]$ is precisely the condition that $\lambda(x)>0$ from Theorem \ref{gammaclassificationbottom}. This shows that $x$ is a $C^1$ function of $(u,v)$.

From the formulas \eqref{xpartials} and the chain rule, we get
$$ \flag_u = \big(\boldsymbol{\eta}'(x) - v \darboux'(x)\big) x_u = \frac{\boldsymbol{\eta}'(x) - v\big(\binormal'(x)-\alpha'(x)\boldsymbol{\eta}'(x)-\alpha(x)\boldsymbol{\eta}''(x)\big)}{1+v\alpha'(x)} = \boldsymbol{\eta}'(x), $$
using the third Frenet-Serret equation \eqref{frenetserreteqs} to eliminate the derivative of $\binormal$.
Similarly we compute
$$ \flag_v = \big(\boldsymbol{\eta}'(x)-v\darboux'(x)\big) x_v - \darboux(x) = -\alpha(x) \flag_u + \alpha(x) \boldsymbol{\eta}'(x) - \binormal(x) = -\binormal(x),$$
using the fact that $x_v = -\alpha x_u$ and the definition of $\darboux$. Since $\boldsymbol{\eta}'=\tangent$ and $\binormal$ are orthonormal at every $x$, we see that $\flag_u$ and $\flag_v$ are orthonormal at every $(u,v)$ in the unit square.

Since $\flag_u=\boldsymbol{\eta}'(x)$ is a composition of the $C^1$ function $\boldsymbol{\eta}'$ and the $C^1$ function $x$, it is also $C^1$. Similarly $\flag_v = -\binormal(x)$ is $C^1$, and this implies that $\flag$ is $C^2$.

Finally we verify the boundary conditions. Since $\alpha(0)=0$, we have $x=0$ whenever $u=0$ in equation \eqref{xcondition}, so that
$$\flag(0,v) = \boldsymbol{\eta}(0) + v\alpha(0) \boldsymbol{\eta}'(0) - v \binormal(0) = v\hat{\jmath},$$
since we constructed $\boldsymbol{\eta}$ to ensure $\binormal(0)=-\hat{\jmath}$. Because $\flag_u(0,v) = \boldsymbol{\eta}'(0) = \hat{\imath}$ for all $v$, the flag is indeed fastened in the horizontal direction all along the flagpole.
\end{proof}

\begin{remark}\label{whipremark1}
In the special case where $\alpha\equiv 0$, the bottom curve $\boldsymbol{\eta}$ remains planar, in the $\hat{\imath}$-$\hat{k}$ plane. In the spherical coordinates of Lemma \ref{frenetserretlemma},
we have $\phi\equiv \beta\equiv 0$, with $\theta'(x) = \kappa(x)$ and $\theta(0)=0$ determining the curve completely. Hence the binormal is constant and given by $\binormal(x) = -\hat{\jmath}$. In addition $x(u,v)$ determined by \eqref{xcondition} is given simply by $x(u,v) = u$.
Hence the parameterization \eqref{regupflagformula} becomes $\flag(u,v) = \boldsymbol{\eta}(u) + v \hat{\jmath}$. In other words, the planar curve at the bottom is vertically translated to fill out the flag.
\end{remark}

\subsection{The manifold structure of the space of upturned flags}
By Theorem~\ref{flagsolutionthm} a regular upturned flag is completely determined by the $C^0$ function $\kappa\colon [0,1]\to\mathbb{R}$ and the $C^1$ function $\alpha\colon [0,1]\to \mathbb{R}$ satisfying the conditions
\begin{equation}\label{alphaconditionsagain}
\alpha(0)=0, \quad \alpha(1)>0, \quad \max\{1-x, \alpha(x)\} + (1-x)\alpha'(x) > 0 \; \forall x\in [0,1].
\end{equation}
Define our linear space containing the $\alpha$ functions to be the space
\begin{equation}\label{alphaspacedef}
\spacey =  \big\{ \alpha\in C^1([0,1]) \,\big\vert\, \alpha(0)=0\big\}, \text{ with Banach norm }
\lVert \alpha\rVert_{\spacey} = \sup_{x\in [0,1]} \lvert \alpha'(x)\rvert,
\end{equation}
We will prove that the conditions in \eqref{alphaconditionsagain} describe an open subset of this space, thereby obtaining the following result concerning the manifold structure of the space of upturned flags.
\begin{theorem}[Manifold structure of regular, upturned flags]\label{flag_manifold}
The space  $\opey$ of functions satisfying the condition \eqref{alphaconditionsagain},
\begin{equation}\label{opensetdef}
\opey = \Big\{\alpha \in \spacey \, \vert \, \alpha(1)>0, \quad \max\{1-x, \alpha(x)\} + (1-x)\alpha'(x) > 0 \; \forall x\in [0,1].\Big\},
\end{equation}
is an open subset of $\spacey$ given by \eqref{alphaspacedef} and therefore a Banach manifold modeled on $\spacey$.

Furthermore, the space $\flagsy$ of regular upturned flags is  diffeomorphic to $C([0,1])\times \opey$, which is an open subset of the Banach space $C([0,1])\times \spacey$ and thus a Banach manifold.
\end{theorem}
\begin{proof}
We have shown in Theorem~\ref{flagsolutionthm} that  regular, upturned flags are uniquely determined by the functions $\kappa\in C([0,1])$ and $\alpha\in  \opey$ and thereby we have established the identification with the set $C([0,1])\times \opey$. It remains to show that $\opey$ is an open subset of  $\spacey$. Therefore, let $\alpha\in \opey$. By Theorem \ref{gammaclassificationbottom}, there is a unique point $x^*\in (0,1)$ such that $\alpha(x)< 1-x$ for $x< x^*$ and $\alpha(x)> 1-x$ for $x>x^*$. Consider a function $f\in \spacey$; we will show that if $\lVert f\rVert_{\spacey}$ is sufficiently small, then $\alpha+f\in \opey$.

For $x\in [0,x^*]$ we have $1+\alpha'(x)>0$, and in particular there is an $\varepsilon_1>0$ such that $1+\alpha'(x)\ge \varepsilon_1$ for $x\in [0,x^*]$. Thus we have
$$ 1+(\alpha + f)'(x) \ge 1+\alpha'(x) - \lVert f\rVert_{\spacey} \ge \varepsilon_1-\lVert f\rVert_{\spacey},$$
which remains positive on $[0,x^*]$ as long as $\lVert f\rVert_{\spacey} < \varepsilon_1$.

For $x\in [x^*,1]$ we similarly have $\alpha(x) + (1-x) \alpha'(x)\ge \varepsilon_2$ for some $\varepsilon_2>0$, and thus
\begin{align*}
(\alpha+f)(x) + (1-x)(\alpha+f)'(x) &\ge \alpha(x) + (1-x) \alpha'(x) - \sup_{x\in [x^*,1]} \lvert f(x)\rvert - \sup_{x\in [x^*,1]} (1-x) \lvert f'(x)\rvert\\
&\ge \varepsilon_2 - 2\lVert f\rVert_{\spacey},
\end{align*}
which remains positive as long as $\lVert f\rVert_{\spacey}< \varepsilon_2/2$.

Requiring that $\lVert f\rVert < \min\{\varepsilon_1, \varepsilon_2/2\}$ ensures that on either interval $[0,x^*]$ or $[x^*,1]$, at least one of the functions is positive, and thus their maximum is also positive. Thus $\alpha+f\in \opey$.
\end{proof}


\begin{remark}\label{whipremark2}
The special case where $\alpha\equiv 0$ generates two-dimensional whips, as mentioned in Remark \ref{whipremark1}. This is obviously a submanifold of $\flagsy$, with tangent space consisting of arbitrary functions $\dot{\kappa}$ with $\dot{\alpha}\equiv 0$. We will show later that in the kinetic energy metric induced on flags, the space of whips is totally geodesic.
\end{remark}

\subsection{The space of smooth, regular upturned flags as a submanifold}\label{submanifoldsection}
Next we would like to consider the space of regular upturned flags as a submanifold of the space of all $C^2$ surfaces. Unfortunately this does not seem to work,
in part due to the complicated smoothness conditions on flags themselves (a $C^0$ function $\kappa$ and a $C^1$ function $\alpha$ generate a $C^2$
curve $\boldsymbol{\eta}$, but not every $C^2$ curve $\boldsymbol{\eta}$ automatically has a $C^1$ function $\alpha$, and the smoothness of the flag surface $\flag$ is even
more involved). Even without these difficulties, the fundamental problem is the same one that arose in Nash's study of isometric immersions~\cite{nash1956imbedding}: the loss of derivatives in the function that maps a parameterized surface to
the induced metric. Here we would like to say that the metric map that takes a parameterized surface $\flag\colon [0,1]^2\to\mathbb R^3$ to its Riemannian metric coefficients via
\begin{equation}\label{isometrymap}
\mathcal{I}\colon C^{k+1}([0,1]^2,\mathbb R^3)\to C^k([0,1]^2, \mathbb R^3), \qquad
\mathcal{I}(\flag) = \big( \tfrac{1}{2} \flag_u\cdot \flag_u,  \tfrac{1}{2} \flag_v\cdot \flag_v, \flag_u\cdot \flag_v\big)
\end{equation}
has $(\tfrac{1}{2}, \tfrac{1}{2},0)$ as a regular value. This requires that the derivative of $\mathcal{I}$ be surjective for any flag, and in particular for the regular upturned ones. In the proposition below we will see that this works in the smooth category ($k=\infty$), but not for flags of finite regularity.
\begin{proposition}\label{derivativesubmersion}
The differential of the map $\mathcal{I}$ defined in \eqref{isometrymap}, at a regular upturned flag $\flag$ as in Definition \ref{regupflagdef} and parameterized as in Theorem \ref{flagsolutionthm}, is given in $(x,y)$ coordinates by
\begin{equation}\label{isometryderivative}
\begin{split}
D\mathcal{I}_{\flag}(z) &= \Big( \frac{f_x(x,y)-\kappa(x) g(x,y)}{1+y\alpha'(x)},
-h_y(x,y) + \alpha(x) \, \frac{h_x(x,y)-\kappa(x)\alpha(x)g(x,y)}{1+y\alpha'(x)}, \\
&\qquad\qquad f_y(x,y) - \alpha(x)\,\frac{f_x(x,y)-\kappa(x)g(x,y)}{1+y\alpha'(x)} - \frac{h_x(x,y)-\kappa(x)\alpha(x)g(x,y)}{1+y\alpha'(x)}\Big),\\
&\qquad\text{ where } z(x,y) = f(x,y) \tangent(x) + g(x,y) \normal(x) + h(x,y) \binormal(x).
\end{split}
\end{equation}
It has a right-inverse given for functions $(p,r,q) = D\mathcal{I}_{\flag}(f,g,h)$ by
\begin{align}
f(x,y) &= \int_0^y \big(2\alpha(x) p(x,s) + q(x,s)\big)\,ds \\
&\qquad\qquad+ \int_0^y \frac{s-y}{1+s\alpha'(x)} \Big( \alpha(x)^2 p_x(x,s)+ \alpha(x) q_x(x,s) + r_x(x,s)\Big) \, ds \label{feq} \\
g(x,y) &= \frac{1}{\kappa(x)} \Big( f_x(x,y) - \big(1+y\alpha'(x)\big) p(x,y) \Big) \label{geq}\\
h(x,y) &= \alpha(x) f(x,y) - \int_0^y \Big( \alpha(x)^2 p(x,s) + \alpha(x) q(x,s) + r(x,s)\Big) \, ds,\label{heq}
\end{align}
which exists for every $(p,q,r)$ iff $\kappa$ is nowhere zero.
\end{proposition}

\begin{proof}[Proof of Proposition~\ref{derivativesubmersion}]
The derivative of $\mathcal{I}$ is given by
\begin{equation}\label{Iderivative}
\begin{split}
D\mathcal{I}_{\flag}(z) = \frac{\partial }{\partial t}\Big|_{t=0} \mathcal{I}\big(\flag(t),\flag(t)\big)
&= \big(\flag_u\cdot \flag_{tu}, \flag_v\cdot \flag_{tv},  \flag_u\cdot \flag_{tv}
 + \flag_{tu}\cdot \flag_v \big)\Big|_{t=0} \\
&= \big( \flag_u\cdot z_u, \flag_v\cdot z_v, \flag_u\cdot z_v + \flag_v\cdot z_u, \big),
\end{split}
\end{equation}
where $z = \frac{\partial \flag}{\partial t}\big|_{t=0}$.

Using the chain rule formulas
$$ \frac{\partial}{\partial u} = \frac{1}{1+y\alpha'(x)} \, \frac{\partial}{\partial x}, \qquad \frac{\partial}{\partial v} = \frac{\partial}{\partial y} - \frac{\alpha(x)}{1+y\alpha'(x)} \, \frac{\partial}{\partial x}$$
as in equation \eqref{xpartials}, along with the fact from Theorem \ref{flagsolutionthm} that
$\flag_u(u,v) = \tangent(x)$ and $\flag_v(u,v) = -\binormal(x)$, the equation \eqref{Iderivative} has components
\begin{multline}\label{isometryderivxy}
D\mathcal{I}_{\flag}(z) = \bigg( \frac{\tangent(x)\cdot z_x(x,y)}{1+y\alpha'(x)},
-\binormal(x)\cdot \Big( z_y(x,y) - \frac{\alpha(x) z_x(x,y)}{1+y\alpha'(x)}, \\
\tangent(x)\cdot \Big( z_y(x,y) - \frac{\alpha(x) z_x(x,y)}{1+y\alpha'(x)}\Big) - \binormal(x) \cdot \frac{z_x(x,y)}{1+y\alpha'(x)}
\Big)\bigg)
\end{multline}
Writing $z$ in the Frenet-Serret basis as in \eqref{isometryderivative} and using \eqref{frenetserreteqs}, we obtain
\begin{align*}
z_x &= (f_x - \kappa g) \tangent + (g_x + \kappa f + \kappa \alpha h) \normal + (h_x - \kappa \alpha g) \binormal \\
z_y &= f_y \tangent + g_y \normal + h_y \binormal,
\end{align*}
and plugging these into \eqref{isometryderivxy} yields the equation \eqref{isometryderivative} for the derivative.

Using formula \eqref{isometryderivative}, we find the right-inverse operator by solving the system
\begin{align}
\frac{f_x(x,y)-\kappa(x) g(x,y)}{1+y\alpha'(x)} &= p(x,y) \label{peq} \\
-h_y(x,y) + \alpha(x) \, \frac{h_x(x,y)-\kappa(x)\alpha(x)g(x,y)}{1+y\alpha'(x)} &= r(x,y)\label{req} \\
f_y(x,y) - \alpha(x)\,\frac{f_x(x,y)-\kappa(x)g(x,y)}{1+y\alpha'(x)} - \frac{h_x(x,y)-\kappa(x)\alpha(x)g(x,y)}{1+y\alpha'(x)} &= q(x,y) \label{qeq}
\end{align}
Multiplying \eqref{qeq} by $\alpha(x)$ and using equations \eqref{peq} and \eqref{req} to replace the $x$-derivatives in it, we obtain
$$ \alpha f_y - \alpha^2 p - (h_y + r) = \alpha q,$$
which is equivalent to
\begin{equation}\label{partialyfheq}
\frac{\partial}{\partial y}\Big( h(x,y) - \alpha(x) f(x,y)\Big) = -F(x,y), \qquad F(x,y) = \alpha(x)^2 p(x,y) + \alpha(x) q(x,y) + r(x,y).
\end{equation}
We can solve \eqref{partialyfheq} for $h$ to get \eqref{heq}, up to an arbitrary function of $x$ which we set to zero.

Then differentiating \eqref{heq} with respect to $x$ and using \eqref{peq} gives
$$ h_x(x,y) - \kappa(x)\alpha(x)g(x,y) = \alpha'(x)f(x,y) + \alpha(x)p(x,y) \big(1+y\alpha'(x)\big) - \int_0^y F_x(x,s)\,ds,$$
and inserting this into \eqref{qeq} gives a single equation for $f$ alone:
$$ f_y(x,y) - \frac{\alpha'(x)}{1+y\alpha'(x)} f(x,y) = 2\alpha(x)p(x,y) + q(x,y) - \frac{1}{1+y\alpha'(x)} \int_0^y  F_x(x,s)\,ds.$$
This is an ordinary differential equation in $y$, which can be solved assuming $f(x,0)=0$ as
$$ f(x,y) = \big( 1 + y\alpha'(x)\big) \int_0^y \frac{2\alpha(x)p(x,s)+q(x,s)}{1+s\alpha'(x)} \, ds + \int_0^y F_x(x,s)\,ds
- \frac{1+y\alpha'(x)}{\alpha'(x)} \int_0^y \frac{F_x(x,s)}{1+s\alpha'(x)} \, ds.$$
Straightforward manipulations using the definition of $F$ in \eqref{partialyfheq} turn this into \eqref{feq}.

Having obtained $f$, we know $h$ from \eqref{heq}. If $\kappa$ is nonzero, we can solve equation \eqref{peq} for $g$.
\end{proof}

\begin{remark}[Loss of derivative]
Note the loss of derivatives in the formulas \eqref{feq}--\eqref{heq}. If we want to show that $D\mathcal{I}_{\flag}$ is surjective
from $C^{k+1}$ surfaces to $C^k$ metric components, then given any $C^k$ functions $(p,r,q)$, we want the solution $(f,g,h)$ to be $C^{k+1}$.
However the formula \eqref{feq} shows that in fact $f$ is only $C^{k-1}$, while $h$ is also $C^{k-1}$ and $g$ is $C^{k-2}$. This observation prevents us from using the inverse function theorem for Banach spaces to show that $C^k$ flags  form a smooth submanifold of $C^k$ surfaces. 
Next we show that this difficulty can be overcome in the smooth category, i.e., we will use the Nash-Moser inverse function theorem to prove that $C^{\infty}$ flags for which the curvature $\kappa$ is nowhere zero form a smooth submanifold of the space of $C^{\infty}$ surfaces. These results are in accordance with the space of volume preserving diffeomorphisms as a sub-Lie group of the full diffeomorphism group, and with the space of regular, volume preserving embeddings as a submanifold of all  regular embeddings~\cite{molitor2017remarks,bauer2016riemannian}.
\end{remark}

\begin{theorem}[Submanifold structure for smooth regular flags]\label{thm:submanifold}
The space of smooth, regular flags with non-vanishing curvature function $\kappa=e(x,0)\neq 0$ is a tame Fr\'echet submanifold
of the space
\begin{equation}
\operatorname{Imm}([0,1]^2,\mathbb R^3)^{\star}:=\left\{\flag \in\operatorname{Imm}([0,1]^2,\mathbb R^3): e(x,0)\neq 0 \right\},
\end{equation}
which is an open subset of the space of all smooth immersions $\operatorname{Imm}([0,1]^2,\mathbb R^3)$.
\end{theorem}
\begin{proof}
Using the results of {Proposition}~\ref{derivativesubmersion}, the proof of this result will follow similarly as in~\cite{molitor2017remarks,bauer2016riemannian} and we will be rather brief in our arguments.
Indeed the situation studied here is much simpler, as the space $C^{\infty}([0,1]^2, \mathbb R^3)$ is a tame Fr\'echet space; in ~\cite{molitor2017remarks,bauer2016riemannian}
the authors consider immersions from a general finite dimensional manifold $M$, which makes the presentation significantly more complicated as it requires one to work in local coordinate charts.

We consider the map $\mathcal I$ in the smooth category:
\begin{equation}
\operatorname{Imm}([0,1]^2,\mathbb R^3)^{\star} \to C^{\infty}([0,1]^2, \mathbb R^3), \qquad
\mathcal{I}(\flag) = \big( \tfrac{1}{2} \flag_u\cdot \flag_u,  \tfrac{1}{2} \flag_v\cdot \flag_v, \flag_u\cdot \flag_v\big).
\end{equation}
We first note that $C^{\infty}([0,1]^2, \mathbb R^3)$ is a tame Fr\'echet space, and that $\operatorname{Imm}([0,1]^2,\mathbb R^3)^{\star} $ is an open subset of it.
Thus in order to apply the Nash-Moser inverse function theorem, it remains to show that there exists an open subset $U\subset \operatorname{Imm}([0,1]^2,\mathbb R^3)^{\star}$ such that
\begin{enumerate}
\item $\mathcal I$ is a smooth, tame map;\label{prop1}
\item $d\mathcal I(x)$ is a linear isomorphism for all $x\in U$;\label{prop2}
\item the map $d\mathcal I^{-1}: U\times C^{\infty}([0,1]^2, \mathbb R^3)\to C^{\infty}([0,1]^2, \mathbb R^3)$ is a smooth tame map. \label{prop3}
\end{enumerate}
Since every nonlinear differential operator is a smooth tame map, see e.g., \cite[Cor. 2.2.7]{hamilton1982inverse}, it follows directly from the definition of $\mathcal I$ that Property~\eqref{prop1} holds. The remaining properties follow directly from the explicit formula for  the inverse $d\mathcal I^{-1}$ given in~Proposition~\ref{derivativesubmersion}. Using that $\operatorname{Imm}([0,1]^2,\mathbb R^3)^{\star}$
is implicitly characterized by the condition $\mathcal I(\flag)=(1/2,1/2,0)$, the result then follows by the Nash-Moser version of the regular value theorem.
\end{proof}

\subsection{Downturned and balanced flags}\label{sec:downturned}
Using the exact same methods, the analogous results also hold for the spaces of downward and balanced flags. The main difference can be seen in the following result, which is the analogue of Theorem~\ref{gammaclassificationbottom}:
\begin{corollary}\label{gammaclassificationtop}
A function $\alpha_t\colon [0,1]\to\mathbb{R}$ with $\alpha_t(0)=0$ and $\alpha_t(1)<0$ generates a family of nonintersecting asymptotic curves in the downturned flag case with diffeomorphic coordinate transformation \eqref{topheavycoords} iff it satisfies the conditions:
 \begin{equation}\label{lambdapositivetop}
 \lambda_t(x)>0  \text{ for all $x\in[0,1]$,}
 \end{equation}
 where
\begin{equation}\label{lambdagammatopdef}
\lambda_t(x) = 1+\alpha_t'(x)\gamma_t(x), \qquad
\gamma_t(x) = \begin{cases} 1 & x-\alpha_t(x)\le 1, \\
-\frac{(1-x)}{\alpha_t(x)} & x-\alpha_t(x)>1.\end{cases}
\end{equation}

In the balanced case, a function $\alpha_b\colon[0,1]\to\mathbb{R}$ with $\alpha_b(0)=\alpha_b(1)=0$ generates nonintersecting asymptotic curves with either \eqref{bottomheavycoords} generating a diffeomorphism on the square if and only if $\lambda_b(x) := 1+\alpha_b'(x)>0$ for all $x\in [0,1]$; here $\gamma_b(x)\equiv 1$ for all $x\in [0,1]$. Equivalently $\alpha_t\colon [0,1]\to \mathbb{R}$ with $\alpha_t(0)=\alpha_t(1)=0$ generates a diffeomorphism via \eqref{topheavycoords} if and only if $\lambda_t(x) := 1-\alpha_t'(x)>0$ for all $x\in [0,1]$.
\end{corollary}
This naturally leads to the following definition.
\begin{definition}\label{regdownflagdef}
A \emph{regular downturned flag} is a $C^2$ isometric embedding $\flag\colon [0,1]^2 \to \mathbb{R}^3$ satisfying the conditions
\begin{equation}\label{flagboundaryconditionsdown}
\flag(0,v)=v \hat{\jmath}, \quad \flag_u(0,v) = \hat{\imath}, 
\quad \text{for all $v\in [0,1]$,}
\end{equation}
and such that there is a $C^1$ function $\alpha_t\colon [0,1]\to \mathbb{R}$ satisfying
$\flag_{uv}(x,1) = -\alpha_t(x) \flag_{uu}(x,1)$ for all $x\in[0,1]$ as well as the conditions of Theorem \ref{gammaclassificationtop}. That is, $\alpha_t(0)=0$, $\alpha_t(1)<0$, and
\begin{equation}\label{nondegenerateflagconddown}
\max\{1-x, -\alpha_t(x)\} + (1-x)\alpha_t'(x) > 0 \quad \text{for all $x\in [0,1]$.}
\end{equation}
\end{definition}

For a regular balanced flag, we can use either this definition or the definition for upturned flags, with the only modifications being that $\alpha_b(1)=\alpha_t(1)=0$ and $1+\alpha_b'(x)>0$ for all $x \in [0,1]$. From the above definition and corollary, it is clear that everything we did for upturned flags can be done in a similar way for both downturned and  balanced flags. For the latter case, the analysis will be significantly easier.

\subsection{The tangent space}\label{tangentspacesection}
\begin{figure}[!ht]
\begin{center}
\includegraphics[width=.45\textwidth]{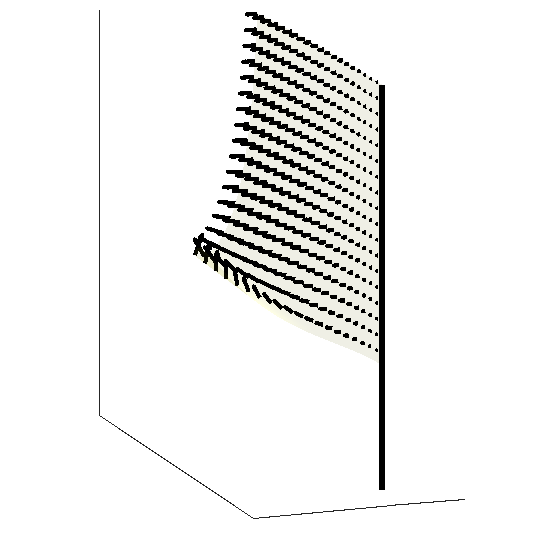}
\includegraphics[width=.45\textwidth]{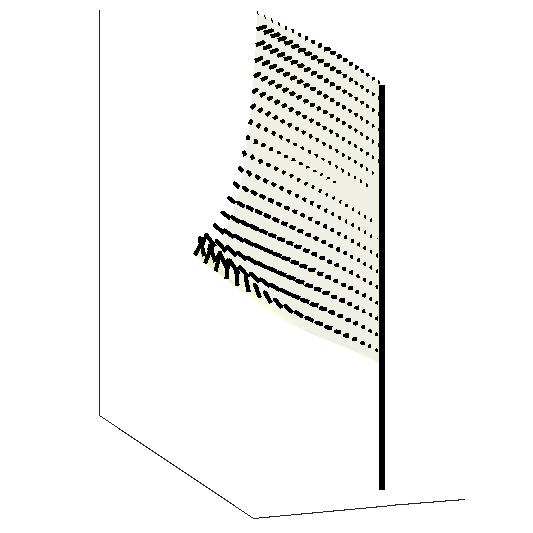}
\caption{Two examples of tangent vectors (vector fields) to the space of regular, upturned flags.}\label{fig:tangentspace}
\end{center}
\end{figure}

From here on we will continue to work again in the finite regularity regime and disregard the submanifold result from Section \ref{submanifoldsection}.
We now compute tangent vectors to the space of flags by considering a curve in the space of flags $\flag(t,u,v)$, and differentiating with respect to $t$.
By Theorem \ref{flagsolutionthm}, this corresponds to a pair of time-dependent functions $\kappa(t,x)$ and $\alpha(t,x)$, which generate a time-dependent bottom edge $\boldsymbol{\eta}(t,x)$ through the coordinates $\theta(t,x)$, $\phi(t,x)$, and $\beta(t,x)$ satisfying the spatial equations \eqref{thetaeq}--\eqref{betaeq} for each fixed $t$. As this notation gets somewhat complicated, we will consider variations using the dot notation: e.g.,
\begin{equation}\label{dotnotation}
\dot{\kappa}(x) = \frac{\partial}{\partial t} \tilde{\kappa}(t,x)\Big|_{t=0}, \qquad \text{where $\tilde{\kappa}(0,x) = \kappa(x)$}.
\end{equation}
In other words, to compute the tangent space at a given regular upturned flag generated by functions $\kappa(x)$ and $\alpha(x)$, we extend to curves $\tilde{\kappa}(t,x)$ and $\tilde{\alpha}(t,x)$ in the space of functions passing through the functions at $t=0$, and find equations for their velocities at time $t=0$.
An example, using the representation of Proposition~\ref{tangentspaceprop}, can be seen in Figure~\ref{fig:tangentspace}. Note that we will use subscript notation for derivatives, and the reader should not confuse the time derivative $\alpha_t(t,x)$ of a bottom-edge $\alpha$ with $\alpha_t(x)$, the top-edge $\alpha$ for a downturned flag. Here all flags are upturned, and from now on it will only represent what we called $\alpha_b$ earlier. 

\begin{proposition}\label{tangentspaceprop}
Suppose $\kappa, \alpha\colon [0,1]\to \mathbb{R}$ are $C^0$ and $C^1$ functions respectively, satisfying the conditions of Definition \ref{regupflagdef}
to generate a regular upturned flag through Lemma \ref{frenetserretlemma} and Theorem \ref{flagsolutionthm}. Let $\dot{\kappa}$ and $\dot{\alpha}$ be $C^0$ and $C^1$ variations, with corresponding variations $\dot{\theta}$, $\dot{\phi}$, $\dot{\beta}$ of the functions in Lemma \ref{frenetserretlemma}. Then the tangent vector to the flag is given by
\begin{equation}\label{tangentflag}
\dot{\flag}(u,v) = \int_0^x \big( \omega(s) \normal(s) + \psi(s) \binormal(s)\big)\,ds + v \Big(
\frac{\psi'(x)}{\kappa(x)} \normal(x) + \psi(x)\big(\tangent(x) + \alpha(x)\binormal(x)\big)\Big),
\end{equation}
where $x=x(u,v)$ is the function solving \eqref{xcondition}.
Here the functions $\omega$ and $\psi$ are related to the variations $\dot{\kappa}$ and $\dot{\alpha}$ by solving the ODEs
\begin{alignat}{3}
\omega'(x) &= \dot{\kappa}(x) - \kappa(x) \alpha(x) \psi(x) &\qquad \omega(0) &= 0  \label{omegaeq} \\
\psi'(x) &= \kappa(x)\chi(x) + \kappa(x) \alpha(x) \omega(x) & \chi(0) &= 0 \label{chieq} \\
\chi'(x) &= -\kappa(x)\dot{\alpha}(x) - \alpha(x) \dot{\kappa}(x) - \kappa(x)\psi(x) & \psi(0)&= 0 \label{psieq}
\end{alignat}
In particular if $\kappa$ and $\dot{\kappa}$ are $C^0$ and $\alpha$ and $\dot{\alpha}$ are $C^1$, then $\omega$, $\psi$, and $\chi$ are all $C^1$.
\end{proposition}

\begin{proof}
By formulas \eqref{regupflagformula} and \eqref{xcondition} in Theorem \ref{flagsolutionthm}, we can write
\begin{equation}\label{flagtimedep}
\flag(t,u,v) = \boldsymbol{\eta}(t,x(t,u,v)) + v \big( \alpha(t,x) \boldsymbol{\eta}_x(t,x) - \binormal(t,x)\big),
\end{equation}
where $x(t,u,v)$ is defined to be the unique solution in $[0,1]$ of
\begin{equation}\label{xtimedep}
x(t,u,v) + \alpha\big(t,x(t,u,v)\big) v = u.
\end{equation}
Differentiating \eqref{flagtimedep} once with respect to $t$ gives (omitting the dependent variables on the right side for brevity):
\begin{equation}\label{rtpart1}
\begin{split}
\flag_t(t,u,v) &= \boldsymbol{\eta}_t + v \big( \alpha_t \boldsymbol{\eta}_x + \alpha \boldsymbol{\eta}_{tx} - \binormal_t) + \Big( \boldsymbol{\eta}_x + v\big( \alpha_x \boldsymbol{\eta}_x + \alpha \boldsymbol{\eta}_{xx} - \binormal_x\big)\Big) \frac{\partial x}{\partial t} \\
&= \boldsymbol{\eta}_t + v \big( \alpha_t \boldsymbol{\eta}_x + \alpha \boldsymbol{\eta}_{tx} - \binormal_t) + (1+v\alpha_x) \boldsymbol{\eta}_x  \frac{\partial x}{\partial t},
\end{split}
\end{equation}
the simplification in the second line being due to the Frenet-Serret equation \eqref{frenetserreteqs}. To find $\frac{\partial x}{\partial t}$, we differentiate \eqref{xtimedep} with respect to $t$ and solve to obtain
\begin{equation}\label{xpartialt}
\frac{\partial x}{\partial t}(t,u,v) = -\frac{v \alpha_t(t,x)}{1+v\alpha_x(t,x)}.
\end{equation}
Using \eqref{xpartialt} in formula \eqref{rtpart1} and simplifying yields
\begin{equation}\label{rtpart2}
\flag_t(t,u,v) = \boldsymbol{\eta}_t(t,x(t,u,v)) + v\, \big( \alpha(t,x(t,u,v)) \boldsymbol{\eta}_{tx}(t,x(t,u,v)) - \binormal_t(t,x(t,u,v)) \big).
\end{equation}
It remains to compute $\binormal_t$ more explicitly.

The formulas \eqref{tangentflag} and \eqref{omegaeq}--\eqref{chieq} are all intrinsic, and can be derived directly from variations
of the Frenet-Serret equations \eqref{frenetserreteqs}. However we will derive them as a consequence of the variations of the coordinate
equations for $\theta$, $\phi$, and $\beta$ given in \eqref{thetaeq}--\eqref{betaeq}, since these are convenient for explicitly constructing the flag numerically.

Differentiating \eqref{thetaeq}--\eqref{betaeq} with respect to time gives
\begin{align}
\dot{\theta}' &= \dot{\kappa} \cos{\beta} - \kappa \dot{\beta} \sin{\beta} \label{dotthetaprimeeq} \\
\dot{\phi}'\cos{\theta} &= \dot{\kappa} \sin{\beta} + \kappa \dot{\beta} \cos{\beta} + \kappa \dot{\theta} \tan{\theta} \sin{\beta}  \label{dotphiprimeeq}\\
\dot{\beta}' &= \dot{\kappa} (\alpha+\sin{\beta} \tan{\theta}) + \kappa (\dot{\alpha} + \dot{\beta} \cos{\beta}\tan{\theta} + \dot{\theta} \sin{\beta} \sec^2{\theta}).\label{dotbetaprimeeq}
\end{align}
Since $\theta$, $\phi$, and $\beta$ are all zero when $x=0$ regardless of time, we must have $\dot{\theta}$, $\dot{\phi}$, and $\dot{\theta}$ also equal to zero when $x=0$.

With $\boldsymbol{\eta}'$ given in terms of $\theta$ and $\phi$ by \eqref{etae1def}, differentiating
 with respect to time gives, in the $\{e_2,e_3\}$ basis of \eqref{e2e3def}, the formula
$$ \dot{\boldsymbol{\eta}}' = -\dot{\theta} e_2 + \dot{\phi} \cos{\theta} e_3,$$
and in terms of the Frenet-Serret basis, we can write this using \eqref{normalbinormaldef} as
\begin{align}
\dot{\boldsymbol{\eta}}'(x) &= \omega(x) \normal(x) + \psi(x) \binormal(x), \label{omegapsidef} \\
\text{where } \dot{\theta}(x) &= \omega(x) \cos{\beta(x)} + \psi(x)\sin{\beta(x)} \label{thetadot} \\
\text{and } \dot{\phi}(x) \cos{\theta(x)} &= \omega(x)\sin{\beta(x)} - \psi(x)\cos{\beta(x)}.\label{phidot}
\end{align}
Since $\dot{\theta}$ and $\dot{\phi}$ are both zero at $x=0$, we find that $\omega(0)=\psi(0)=0$ as well.

Differentiating \eqref{thetadot}--\eqref{phidot} with respect to $x$, we obtain
\begin{align*}
\dot{\theta}' &= (\omega'+\psi \beta') \cos{\beta} + (\psi' - \omega \beta') \sin{\beta}, \\
\dot{\phi}'\cos{\theta} &=(\omega'+\psi\beta') \sin{\beta} - (\psi'-\omega \beta') \cos{\beta} +  \dot{\phi}\theta' \sin{\theta}.
\end{align*}
Matching with \eqref{dotthetaprimeeq}--\eqref{dotphiprimeeq}, we get the system
\begin{equation}\label{dotkappadotbeta}
\begin{split}
\dot{\kappa} \cos{\beta} - \kappa \dot{\beta} \sin{\beta} &= (\omega'+\psi \beta') \cos{\beta} + (\psi' - \omega \beta') \sin{\beta} \\
\dot{\kappa} \sin{\beta} + \kappa \dot{\beta} \cos{\beta} &= (\omega'+\psi\beta') \sin{\beta} - (\psi'-\omega \beta') \cos{\beta} 
-\kappa \dot{\theta} \tan{\theta} \sin{\beta} +  \dot{\phi}\theta'\sin{\theta}.
\end{split}
\end{equation}
Solving for $\dot{\kappa}$ gives
\begin{align*}
\dot{\kappa} &= \omega' + \psi \beta' - \kappa \dot{\theta} \tan{\theta} \sin^2{\beta} + \dot{\phi} \theta' \sin{\theta} \sin{\beta} \\
             &= \omega' + \psi\kappa \alpha,
\end{align*}
using equation \eqref{thetaeq}, \eqref{betaeq}, \eqref{thetadot}, and \eqref{phidot}. This is \eqref{omegaeq}.

Similarly solving the system for $\dot{\beta}$, we get
\begin{align*}
\kappa \dot{\beta} &= -\psi' + \omega \beta' - \kappa \dot{\theta} \tan{\theta} \cos{\beta}\sin{\beta} + \dot{\phi}\theta' \sin{\theta} \cos{\beta} \\
&= -\psi' + \omega\kappa \alpha + \omega\kappa \sin{\beta} \tan{\theta} - \kappa \psi \cos{\beta} \tan{\theta},
\end{align*}
again using \eqref{thetaeq}, \eqref{betaeq}, \eqref{thetadot}, and \eqref{phidot}.
Defining the auxiliary function $\chi$ by the formula \eqref{chieq}, this becomes
\begin{equation}\label{betadotomegachipsi}
\dot{\beta} = -\chi + \tan{\theta}(\omega \sin{\beta} - \psi \cos{\beta}),
\end{equation}
and the fact that $\chi(0)=0$ follows from the fact that $\dot{\beta}(0)=0$ together with $\theta(0)=0$.

Now differentiating equation \eqref{betadotomegachipsi} with respect to $x$ gives
$$\dot{\beta}' = -\chi'
+ \theta' \sec^2{\theta}  (\omega \sin{\beta} - \psi \cos{\beta})
+ \beta' \tan{\theta} (\omega \cos{\beta} + \psi \sin{\beta})
+ \tan{\theta} (\omega' \sin{\beta} - \psi' \cos{\beta}),
$$
and matching equation \eqref{dotbetaprimeeq} for $\dot{\beta}'$ leaves an equation for $\dot{\alpha}$.
We eliminate $\theta'$ from this using \eqref{thetaeq}, $\beta'$ using \eqref{betaeq}, $\omega'$ using \eqref{omegaeq}, $\dot{\theta}$ using \eqref{thetadot}, and $\dot{\beta}$ using \eqref{betadotomegachipsi}. What remains after the cancellations is equation \eqref{psieq}.

Finally we return to the computation of $\binormal_t(t,x)$. By formula \eqref{normalbinormaldef}, we have
$\binormal(t,x) = -\sin{\beta} e_2 - \cos{\beta} e_3$,
so that using formula \eqref{e2e3def}, we obtain
\begin{equation}\label{dotbinormal1}
\begin{split}
\dot{\binormal} &= \dot{\beta} (-\cos{\beta} e_2 + \sin{\beta} e_3) - \sin{\beta} (\dot{\theta} e_1 + \dot{\phi} \sin{\theta} e_3)
+ \dot{\phi} \cos{\beta} (\cos{\theta} e_1 + \sin{\theta} e_2) \\
&= -\psi \tangent - \chi \normal,
\end{split}
\end{equation}
and inserting this into formula \eqref{rtpart2} gives the result \eqref{tangentflag}.
 \end{proof}

\begin{remark}\label{constraintstangent}
Note that $\dot{\kappa}$ is unconstrained since $\kappa$ is thus far unconstrained, while $\dot{\alpha}$ is unconstrained except that
$\alpha(0)=0$, since the nondegeneracy condition \eqref{nondegenerateflagcond} is an open condition in the $C^1$ topology.
However the functions $\omega$ and $\psi$ are constrained: if we wish to solve the system \eqref{omegaeq}--\eqref{psieq} algebraically for $\dot{\kappa}$
and $\dot{\alpha}$ given $\omega$ and $\psi$, we need to worry about any points where $\kappa$ is equal to zero. First we need to solve \eqref{chieq}
for $\chi$, which is only possible if $\psi'=0$ whenever $\kappa=0$, and then we need to ensure that $(\chi' + \alpha \dot{\kappa})=0$ whenever $\kappa=0$. Furthermore even if we could ensure these conditions, they would not necessarily lead to a $C^1$ function $\dot{\alpha}$, unless we knew higher-order derivative conditions on $\kappa$. Later when needed to derive the geodesic equation, we will work formally, assuming either that $\kappa$ is nowhere zero or that the functions $\omega$ and $\psi$ can be specified somewhat arbitrarily, but for the more rigorous analysis of this as an infinite-dimensional geodesic system, one would need to worry about this.
\end{remark}

It will be convenient later, when deriving the geodesic equation, to specify extra smoothness conditions on the flag at the flagpole. Since the flag is constrained to have $\flag(0,v)$ fixed at $(0,v,0)$ for all $v\in [0,1]$, it is natural to demand that the $C^2$ function $\flag$ extend to an odd function in the $u$ variable over $[-1,1]$, which imposes the additional condition that $\flag_{uu}(0,v)=0$ for all $v$, and this is equivalent to requiring that $\kappa(0)=0$. Since $\alpha$ is a $C^1$ function and we have already assumed that $\alpha(0)=0$, there is no additional condition to impose on it.

\begin{proposition}\label{oddflagsprop}
The space $\flagsyodd$ of odd regular upturned flags is defined to be those flags generated via Theorem \ref{flagsolutionthm}
such that $\kappa$ and $\alpha$ extend to odd functions through $x=0$. Its tangent space consists of $C^0$ functions $\dot{\kappa}$ and $C^1$ functions $\dot{\alpha}$ which extend to odd functions through $x=0$, and $\flagsyodd$ is a submanifold of $ \flagsy$.

For odd functions $\kappa$ and $\alpha$, the function $\theta$ given by \eqref{thetaeq} is even, while
$\phi$ and $\beta$ given by \eqref{phieq}--\eqref{betaeq} are odd. Similarly for odd functions $\dot{\kappa}$ and $\dot{\alpha}$,
the function $\omega$ given by \eqref{omegaeq} is even, while the functions $\psi$ and $\chi$ given by  \eqref{chieq}--\eqref{psieq}
are odd.
\end{proposition}

\begin{proof}
Extend the solutions $\theta$, $\phi$, and $\beta$ of \eqref{thetaeq}--\eqref{betaeq} to the interval $[-1,1]$. Define 
$\tilde{\theta}(x) = \theta(-x)$, $\tilde{\phi}(x) = -\phi(-x)$, and $\tilde{\beta}(x)=-\beta(x)$. Then these new functions satisfy the ODEs
\begin{align*}
\tilde{\theta}'(x) &= -\theta'(-x) = -\kappa(-x) \cos{\beta(-x)} = \kappa(x) \cos{\tilde{\beta}(x)} \\
\tilde{\phi}'(x) &= \phi'(-x) = \kappa(-x) \sec{\theta(-x)} \sin{\beta(-x)} = \kappa(x) \sec{\tilde{\theta}(x)} \sin{\tilde{\beta}(x)} \\
\tilde{\beta}'(x) &= \beta'(-x) = \kappa(-x)\big( \alpha(-x)+\sin{\beta(-x)} \tan{\theta(-x)}\big) = \kappa(x) \big( \alpha(x) + \sin{\tilde{\beta}(x)} \tan{\tilde{\theta}(x)}\big),
\end{align*}
using the assumption that $\kappa$ and $\alpha$ are odd. This is the same system as \eqref{thetaeq}--\eqref{betaeq}. 

Since the initial conditions $\tilde{\theta}(0)=\tilde{\phi}(0)=\tilde{\beta}(0)=0$ do not change,
uniqueness of solutions of ODEs implies that $\theta=\tilde{\theta}$, $\phi=\tilde{\phi}$, and $\beta=\tilde{\beta}$.
Thus $\theta$ is even while $\phi$ and $\beta$ are odd.
The statements about $\omega$, $\psi$, and $\chi$ follow the same way from the system \eqref{omegaeq}--\eqref{psieq}.

The submanifold result for $\flagsyodd$ is obvious since the only additional constraint on the space is that
$\kappa(0)=0$, which is a closed linear subspace of the first component.
\end{proof}

\section{The kinetic energy}\label{sec:kinetic_energy}
In this section we will consider a natural Riemannian metric on the space of regular upturned flags, which is induced by the kinetic energy metric on the space of general surfaces.
\begin{definition}\label{riemannianmetricdef}
If $\flag$ is a regular upturned flag as in Definition \ref{regupflagdef}, and $\dot{\flag}$ is a tangent vector
as in  Proposition \ref{tangentspaceprop}, then the kinetic energy Riemannian metric is defined to be
\begin{equation}\label{riemannianmetric}
\llangle \dot{\flag}, \dot{\flag}\rrangle_{\flag} = \int_0^1 \int_0^1 \lvert \dot{\flag}(u,v)\rvert^2 \, du\,dv.
\end{equation}
{The corresponding kinetic energy Lagrangian is then given by
\begin{equation}\label{actionfunctional}
\mathcal L_{\operatorname{kin}}(\flag) = \frac{1}{2} \int_0^T \left\langle\!\!\left\langle \frac{d\flag}{dt},\frac{d\flag}{dt}\right\rangle\!\!\right\rangle^2_{\flag(t,u,v)} \,dt,
\end{equation}
where $\flag$ is a path of flags subject to endpoint conditions $\flag(0) = \flag_0$ and $\flag(T) = \flag_1$ for two given regular upturned flags $\flag_0$ and~$\flag_1$. }
\end{definition}

{The kinetic energy metric (Lagrangian, resp.)} is naturally expressed in the $(u,v)$ coordinates, but more easily computed in the $(x,y)$ coordinates of
formula \eqref{xytouv}, since all the important functions depend only on the $x$ variable.

\begin{proposition}\label{metric_coordinates}
In terms of the functions $\omega$ and $\psi$ defined in Proposition \ref{tangentspaceprop}, and the functions $\gamma$ and $\lambda$ defined in
Theorem \ref{gammaclassificationbottom} by formula \eqref{lambdagammadef},
the Riemannian metric \eqref{riemannianmetric}  is given by
\begin{multline}\label{riemannianmetricxy}
\llangle  \dot{\flag}, \dot{\flag}\rrangle_{\flag} = \int_0^1 \bigg( \frac{\gamma(x)}{2}\, (1+\lambda(x))\, \lvert \dot{\boldsymbol{\eta}}(x)\rvert^2 \\
  + \frac{\gamma(x)^2}{3} \, (1 + 2\lambda(x))\, \left\langle \dot{\boldsymbol{\eta}}(x), \frac{\psi'(x)}{\kappa(x)} \normal(x) + \psi(x) \darboux(x)\right\rangle \\
    + \frac{\gamma(x)^3}{12} \, ( 1 + 3\lambda(x))\, \left( \frac{\psi'(x)^2}{\kappa(x)^2} + \big(1+\alpha(x)^2\big) \psi(x)^2\right) \bigg) \, dx,
\end{multline}
where
$$ \frac{d}{dx} \dot{\boldsymbol{\eta}}(x) = \omega(x) \normal(x) + \psi(x) \binormal(x), \qquad \dot{\boldsymbol{\eta}}(0)=0.$$
\end{proposition}
\begin{remark}
An alternative formula, circumventing the functions $\omega$ and $\psi$ and involving instead the time derivative of $\binormal$ and the space derivative of $\boldsymbol{\eta}$, is given by
\begin{equation}\label{eq:metric_alt}
\begin{aligned}
\\\llangle  \dot{\flag}, \dot{\flag}\rrangle_{\flag}=&\int \frac{\gamma(x)}{2}\, (1+\lambda(x))\lvert \dot{\boldsymbol{\eta}}(x)\rvert^2\\&\qquad+ (\alpha(x)\langle \dot{\boldsymbol{\eta}}(x),\tfrac{d}{dx}\dot{\boldsymbol{\eta}}(x)\rangle-\langle \dot{\boldsymbol{\eta}}(x),\dot \binormal(x)\rangle) \frac{\gamma(x)^2}{3} \, (1 + 2\lambda(x))
      \\&\qquad+(\alpha^2(x) |\tfrac{d}{dx}\dot{\boldsymbol{\eta}}(x)|^2-2\alpha(x)\langle\tfrac{d}{dx}\dot{\boldsymbol{\eta}}(x) ,\dot\binormal(x)\rangle+ |\dot\binormal(x)|^2) \frac{\gamma(x)^3}{12} \, ( 1 + 3\lambda(x))dx.
\end{aligned}
\end{equation}
\end{remark}
\begin{proof}
Using formula \eqref{tangentflag}, we have that
$$\dot{\flag}(u,v) = \dot{\boldsymbol{\eta}}(x) + y \left(\frac{\psi'(x)}{\kappa(x)} \normal(x) + \psi(x) \dorbaux(x)\right),$$
where $\dot{\boldsymbol{\eta}}'(x) = \omega(x) \normal(x) + \psi(x) \binormal(x)$ and $\dorbaux(x) = \tangent(x) + \alpha(x) \binormal(x)$,
and $(x,y)$ are related to $(u,v)$ by the formula $v=y$, $u=x+\alpha(x)y$. The area forms are related by the Jacobian \eqref{jacobian}:
\begin{equation}\label{jacobianarea}
du\wedge dv = \big(1+y\alpha'(x)\big) \, dx\wedge dy,
\end{equation}
and the right side is always positive for $0\le y\le \gamma(x)$ by the assumption \eqref{lambdapositive}.

Applying the change of variables, we then get
\begin{align*}
\llangle \dot{\flag}, \dot{\flag}\rrangle_{\flag} &=
\int_0^1 \int_0^{\gamma(x)} \Big( \big(1+y\alpha'(x)\big) \lvert \dot{\boldsymbol{\eta}}(x)\rvert^2 +
   2y(1+y\alpha'(x)) \left\langle \dot{\boldsymbol{\eta}}(x), \frac{\psi'(x)}{\kappa(x)} \normal(x) + \psi(x) \dorbaux(x)\right\rangle \\
&\qquad\qquad+ y^2 (1+y\alpha'(x)) \left\lvert \frac{\psi'(x)}{\kappa(x)} \normal(x) + \psi(x) \dorbaux(x)\right\rvert^2 \Big) \, dy \, dx  \\
&= \int_0^1 \Big( \big( \gamma(x) + \tfrac{1}{2} \gamma(x)^2 \alpha'(x)\big)\lvert \dot{\boldsymbol{\eta}}(x)\rvert^2  \\
 &\qquad\qquad+ \big( \gamma(x)^2 + \tfrac{2}{3} \gamma(x)^3 \alpha'(x)\big)  \left\langle \dot{\boldsymbol{\eta}}(x), \frac{\psi'(x)}{\kappa(x)} \normal(x) + \psi(x) \dorbaux(x)\right\rangle \\
 &\qquad\qquad+ \big( \tfrac{1}{3} \gamma(x)^3 + \tfrac{1}{4}  \gamma(x)^4\alpha'(x) \big) \left( \frac{\psi'(x)^2}{\kappa(x)^2} + \big(1+\alpha(x)^2\big) \psi(x)^2\right) \Big) \, dx.
\end{align*}
This can easily be simplified using $\lambda(x) = 1+\alpha'(x) \gamma(x)$ to the formula \eqref{riemannianmetricxy}.
\end{proof}

While obviously quite complicated, the formula \eqref{riemannianmetricxy} has the advantage that it involves only functions of the $x$ variable,
and thus it represents a Riemannian metric directly on the space of unit-speed curves $\boldsymbol{\eta}$. More explicitly, since it is obviously quadratic in the velocity components $\omega$ and $\psi$, and since those depend in a linear (albeit very nonlocal) way on the functions $\dot{\kappa}$ and $\dot{\alpha}$ through the equations \eqref{omegaeq}--\eqref{psieq}, we obtain a highly nonlocal Riemannian metric on the manifold $\flagsy$ defined by Theorem \ref{flag_manifold}.

The reason this is useful is because we may then construct solutions of the boundary-value problem by minimizing the {Lagrangian~\eqref{actionfunctional}}.
Conceptually it is easy to consider an algorithm that chooses intermediate functions $\kappa(t_i)$, $\alpha(t_i)$, $\dot{\kappa}(t_i)$, and $\dot{\alpha}(t_i)$ for a partition $\{t_0, \ldots, t_m\}$ of $[0,T]$ in order to minimize the total action, although the actual computations to do this involve numerically solving the ODEs \eqref{thetaeq}--\eqref{betaeq} and \eqref{omegaeq}--\eqref{psieq} for each fixed time $t_i$ in order to be able to plug in to the action functional \eqref{actionfunctional}. We will follow this approach in Section~\ref{sec:numerics}, where we will present selected numerical experiments.

\begin{remark}\label{whipremark3}
In the special case of whips, as in Remark \ref{whipremark1} and Remark \ref{whipremark2}, we have $\alpha\equiv 0$ and $\dot{\alpha}\equiv 0$.
As a result we get $\gamma(x)\equiv 1$ and $\lambda(x)\equiv 1$, from the definitions \eqref{lambdagammadef}. Furthermore by \eqref{chieq}--\eqref{psieq},
we have that $\psi$ and $\chi$ satisfy the system
$$ \psi'(x) = \kappa(x) \chi(x), \qquad \chi'(x) = -\kappa(x) \psi(x), \qquad \psi(0)=\chi(0)=0,$$
whose unique solution is $\psi\equiv \chi \equiv 0$.
The formula \eqref{riemannianmetricxy} thus simplifies to
\begin{equation}\label{whipriemannian}
\llangle \dot{\flag}, \dot{\flag}\rrangle_{\flag} = \int_0^1 \lvert \dot{\boldsymbol{\eta}}(x)\rvert^2 \, dx,
\end{equation}
which is the usual kinetic energy for the space of two-dimensional inextensible curves. This shows that the space of whips is an isometrically embedded
submanifold of the space of regular upturned flags.
\end{remark}

\subsection{{Including the effects of gravity and wind}\label{sec:gravity}}
{Next we will describe how one could include the {external} effects of gravity {and wind} by including {extra terms in the Lagrangian.}
\begin{definition}\label{riemannianmetricdefpotential}
Let $\flag$ be a regular upturned flag as in Definition \ref{regupflagdef}. 
Then the gravitational energy is defined to be
\begin{equation}\label{Egrav}
\mathcal E_{\operatorname{Grav}}(\flag)=\int_0^1\int_0^1 \flag(u,v) du dv\cdot  \hat{\jmath},
\end{equation}
where, for simplicity, we set the gravitational constant to be equal to one. 
The corresponding gravitational energy Lagrangian is then given by
\begin{equation}\label{grav_energy}
\mathcal L_{\operatorname{Grav}}(\flag) =  \int_0^T \mathcal E_{\operatorname{Grav}}(\flag) dt,
\end{equation}
where $\flag$ is again a path of flags subject to endpoint conditions $\flag(0) = \flag_0$ and $\flag(T) = \flag_1$ for two given regular upturned flags $\flag_0$ and~$\flag_1$. 
\end{definition}}

{
\begin{remark}
Using this definition the motion of a flag considering its kinetic energy and gravity can be described as a solution to the total energy Lagrangian
\begin{equation}
\mathcal L(\flag)=\mathcal L_{\operatorname{kin}}(\flag) -\mathcal L_{\operatorname{Grav}}(\flag), 
\end{equation}
subject to the same boundary conditions as above.
\end{remark}}

{
In the next proposition we calculate an expression for the gravitational energy in the $(x,y)$ coordinates of
formula \eqref{xytouv}:
\begin{proposition}\label{metric_coordinates_potential}
In terms of the curve $\boldsymbol{\eta}$ and the functions $\alpha$ and $\gamma$,
the gravitational energy \eqref{Egrav} of a flag $\flag$ is given by
\begin{equation}\label{Egrav_xy}
\mathcal E_{\text{Grav}}(\flag)=\int_0^1  \boldsymbol{\eta}(x)(\gamma(x)+\tfrac{\gamma(x)^2}2 \alpha'(x))+  \boldsymbol{\eta}'(x)(\tfrac{\gamma(x)^2}2+\tfrac{\gamma(x)^3}3 \alpha'(x))-\binormal(x)(\tfrac{\gamma(x)^2}2+\tfrac{\gamma(x)^3}3 \alpha'(x)) \, dx\cdot  \hat{\jmath}
\end{equation}
\end{proposition}
\begin{proof}
By formula \eqref{regupflagformula} we have that 
$$
\flag(x,y) = \boldsymbol{\eta}(x) + y\alpha(x) \boldsymbol{\eta}'(x)-y\binormal(x).
$$
Using formula~\eqref{jacobianarea} for the Jacobian of the coordinate change, we thus have
\begin{align}
&\mathcal E_{\text{Grav}}(\flag)= \int_0^1\int_0^{\gamma(x)}  (\boldsymbol{\eta}(x) + y\alpha(x) \boldsymbol{\eta}'(x)-y\binormal(x))\big(1+y\alpha'(x)\big) \, dxdy \cdot  \hat{\jmath}\\
&\quad=\int_0^1\int_0^{\gamma(x)}  (\boldsymbol{\eta}(x) + y\alpha(x) \boldsymbol{\eta}'(x)-y\binormal(x))+(y\boldsymbol{\eta}(x) + y^2\alpha(x) \boldsymbol{\eta}'(x)-y^2\binormal(x))\alpha'(x)\big) \, dxdy \cdot  \hat{\jmath}.\end{align}
Now the desired formula follows by integrating in the variable $y$.
\end{proof}}

{
\begin{remark}[Modelling the effects of wind]
The next step to obtain a physically realistic model would be to include the effects of wind, i.e., the interaction of the flag with the surrounding fluid. As compared to the rather simple nature of the gravitational force, this is a much more challenging problem and several approaches have been considered in the literature, see e.g. ~\cite{zhang2000flexible}, \cite{1968JPSJ...24..392T}, \cite{shelley2011flapping}, \cite{argentina2005fluid} or \cite{fitt2001unsteady}. Assuming an ideal fluid with potential flows  one would assume that one could define a new energy term by considering the velocity field generated by the movement of the flag in the fluid and the intrinsic kinetic energy of the fluid. This would, in particular, require one to solve the interface boundary conditions of the flag with the fluid. In future work it would be interesting to also perform a similar analysis for the model of the present article.
 \end{remark}}

\section{The geodesic equation}
In this section we will calculate the geodesic equation of the kinetic energy metric introduced {in the previous section.} These equations can be interpreted as the governing equations for the motion of a flag {(ignoring the effects of gravity and wind). As the resulting formulas are already rather technically involved, we will not present the Euler-Lagrange equation of the total energy, i.e., including gravity and wind. For the case of gravity the derivation would follow similarly; as mentioned previously we believe that adding the effects of wind to this model, while certainly interesting, would be significantly more difficult and is outside the scope of the present article.} 

The geodesic equation is obtained by minimizing the action \eqref{actionfunctional}. We consider a family of regular upturned
flags depending on time and on some small parameter $\varpar$, as $\flag(\varpar, t, u,v)$ for $\varpar\in (-\epsilon, \epsilon)$, $t\in [0,T]$, $u,v\in [0,1]$. Differentiating the action with respect to $\varpar$, we obtain the requirement that
$$ \frac{dS}{d\varpar}\Big|_{\varpar=0} = \int_0^T \left\langle\!\!\left\langle \frac{\partial^2 \flag}{\partial t \partial \varpar}, \frac{\partial \flag}{\partial t}\right\rangle\!\!\right\rangle_{\flag(\varpar, t,u,v)} \,dt = 0
$$
for every variation $\flag(\varpar, t, u,v)$ fixed at the endpoints $t=0$ and $t=T$. Since the Riemannian metric \eqref{riemannianmetric}
does not depend explicitly on the flag $\flag$ when expressed in $(u,v)$ coordinates, we can simply integrate by parts in time to obtain the condition
$$
\int_0^T \int_0^1 \int_0^1 \langle \flag_{tt}(t,u,v), \varfield(t,u,v)\rangle \, du\,dv\, dt = 0, \qquad \varfield(t,u,v) = \frac{\partial \flag(\varpar,t,u,v)}{\partial \varpar}\Big|_{\varpar=0}.
$$
The time-integral formula above is zero for all time-dependent variations $\varfield(t,u,v)$ if and only if the integrand is zero
at each time: that is,
\begin{equation}\label{variationalacceleration}
\int_0^1 \int_0^1 \langle \ddot{\flag}(u,v), \varfield(u,v)\rangle \, du\,dv = 0, \qquad \text{for every variation field $\varfield$.}
\end{equation}
Equation \eqref{variationalacceleration} must hold for every possible choice of $\varfield(u,v)$, which means $\ddot{\flag}$ must be orthogonal to every possible tangent vector at the given flag, all of which are described by Proposition \ref{tangentspaceprop}.

In the following theorem we present these equations for a flag that is either balanced or upturned. In addition we will assume that $\kappa(0)=0$---the oddness condition---to ensure compatibility at $x=0$.
\begin{theorem}[Geodesic equation on the space of odd, upturned flags]\label{balancedgeodesicthm}
Given initial conditions $\flag(0) \in \mathcal F_o$ and $\dot\flag(0)\in T_{\flag(0)}  \mathcal F_o$ that are described by their generating functions $\alpha$ and $\kappa$ ($\dot \alpha$ and $\dot\kappa$ resp.), the geodesic equation on the space of odd, regular, upturned flags is given by the second order equation
\begin{align}
\ddot{\kappa} &= \frac{d^2}{dx^2}(\xi \kappa)  + \frac{d}{dx}\big((z+\alpha \varphi)\kappa\big) + (\alpha q + \chi^2-\psi^2)\kappa \quad \text{and} \label{kappaddot} \\
\ddot{\alpha} &= -\frac{2\dot{\kappa} \dot{\alpha}}{\kappa}
+ \frac{\alpha'}{\kappa} \frac{d}{dx}(\xi\kappa) - \frac{d}{dx} \big( \mu\kappa\big) + \alpha' (z + \alpha \varphi)
- (1+\alpha^2)q + \alpha(\omega^2 - \chi^2) + 2\omega\chi.\label{alphaddot}
\end{align}
where $\omega, \psi, \chi\colon [0,1]\to\mathbb{R}$ are defined by equations \eqref{omegaeq}--\eqref{psieq}, and where the remaining
coefficient functions  $(\sigma, \rho, z, \varphi, g, q)$ are defined as solutions  to the following ODE system with homogeneous boundary conditions on $[0,1]$:
\begin{alignat}{3}
\sigma' &= \frac{\gamma(1+\lambda)}{2} z + \frac{\gamma^2(1+2\lambda)}{6} \Big( q - \alpha(\omega^2 + \psi^2) - 2\omega \chi\Big), & \quad \sigma(1)&=0 \label{balsigma} \\
\rho' &= \frac{\gamma(1+\lambda)}{2} \varphi + \frac{\gamma^2(1+2\lambda)}{6} \Big( \alpha q + \chi^2 + \psi^2\Big), & \rho(1) &= 0 \label{balrho} \\
g' &= -\rho + \frac{\gamma^2(1+2\lambda)}{6} (z + \alpha \varphi) + \frac{\gamma^3(1+3\lambda)}{12} \Big( (1+\alpha^2)q + \alpha (\chi^2-\omega^2) - 2\omega \chi\Big), & g(1)&= 0 \label{balg} \\
z' &= \kappa^2 \xi - (\omega^2 + \psi^2), & z(0) &= 0 \label{balz} \\
\varphi' &= q + \alpha \kappa^2 \xi, & \varphi(0)&=0 \label{balphi} \\
q' &= 2(\chi \dot{\kappa} - \kappa \psi \omega) + \kappa^2 \mu, & q(0)&= 0. \label{balq}
\end{alignat}
Here $\mu$ and $\xi$ are defined via
\begin{alignat}{3}
&\mu = \frac{1}{\gamma^3(\lambda^2 + 4\lambda + 1)} \Big( 36(1+\lambda) g - 12\gamma (1+2\lambda)(\sigma + \alpha \rho) - \gamma^4 (5+\lambda)\dot{\alpha}^2\Big)&& \label{mueq} \\
&\xi = \frac{1}{\gamma^2(\lambda^2+4\lambda+1)} \big( - 12(1+2\lambda) g + 6\gamma (1+3\lambda)(\sigma + \alpha \rho)  + \gamma^4 \dot{\alpha^2}\big).&&\label{xieq}
\end{alignat}
\end{theorem}

Because $\lambda>0$ always by assumption, the equations \eqref{balsigma}--\eqref{xieq} have nonsingular coefficients, except where $\gamma=0$. In the balanced case, this never happens since $\gamma$ is always $1$. In the upturned case the only time $\gamma$ vanishes is at $x=1$, and the equations can be rewritten in a nonsingular way there, which are likely easier to work with numerically.

\begin{proposition}\label{singularat1}
In the case of an upturned flag, the following functions are nonsingular on $[0,1]$:
$$ \tilde{\sigma} = \frac{\sigma}{(1-x)^2}, \qquad \tilde{\rho} = \frac{\rho}{(1-x)^2}, \qquad \tilde{g} = \frac{g}{(1-x)^3},$$
$$ \tilde{\gamma}(x) = \frac{\gamma(x)}{1-x} = \begin{cases} \frac{1}{1-x} & x\le x^* \\
\frac{1}{\alpha(x)} & x>x^*.\end{cases}$$
Equations \eqref{balsigma}--\eqref{xieq} can all be rewritten in terms of them to avoid singular coefficients.
\end{proposition}

\begin{proof} Recall that
$$ \gamma(x) = \begin{cases} 1 & x\le x^*, \\
\frac{1-x}{\alpha(x)} & x>x^*. \end{cases}$$
On the interval $[x^*,1]$, we compute that $\gamma(x)$ is a decreasing function since
$$\gamma'(x) = \frac{d}{dx} \left(\frac{1-x}{\alpha(x)}\right) =  -\lambda(x)/\alpha(x) < 0, \text{ for $x\in (x^*,1]$.}$$
By assumption $\alpha(1)>0$, and so the only time $\gamma$ vanishes is at $x=1$, and when this happens $\gamma'(1)<0$.

Since $\sigma(1)=0$, equation \eqref{balsigma} implies that $\sigma'(1)=0$ as well, and thus
$$ \lim_{x\to 1} \frac{\sigma(x)}{\gamma(x)^2} = \lim_{x\to 1} \frac{\sigma'(x)}{2\gamma(x)\gamma'(x)} =
\frac{1}{2\gamma'(1)} \lim_{x\to 1} \left[ \frac{1+\lambda}{2} z + \frac{\gamma(1+2\lambda)}{6} \Big( q - \alpha(\omega^2 + \psi^2) - 2\omega \chi\Big)\right] = \frac{z(1)}{2\gamma'(1)}.$$
Similarly by \eqref{balrho} we have
$$ \lim_{x\to 1} \frac{\rho(x)}{\gamma(x)^2} = \frac{\varphi(1)}{2\gamma'(1)}.$$
Since $\sigma$ and $\rho$ behave like $(1-x)^2$ near $x=1$, equation \eqref{balg} implies that $g$ behaves like $(1-x)^3$ near $x=1$, and
\begin{align*}
\lim_{x\to 1} \frac{g(x)}{\gamma(x)^3} &= \frac{1}{3\gamma'(1)} \lim_{x\to 1} \left( -\frac{\rho(x)}{\gamma(x)^2} + \frac{1+2\lambda}{6} (z+\alpha\varphi)\right)\\
&= \frac{1}{3\gamma'(1)} \left( \frac{\varphi(1)}{2\gamma'(1)} + \frac{1+2\lambda(1)}{6}\, \big(z(1)+\alpha(1)\varphi(1)\big)\right).
\end{align*}

This means that $\mu$ defined by \eqref{mueq} has a finite limit as $x\to 1$, and that $\xi$ defined by \eqref{xieq} in fact approaches zero as $x\to 1$, which means the terms appearing in \eqref{balz}--\eqref{balq} are continuous on all of $[0,1]$.
\end{proof}

%


Before we  prove Theorem~\ref{balancedgeodesicthm}, we will first review the geodesic equation in the case of whips, which will illustrate the equation in a simpler case and how the oddness assumption of Proposition \ref{oddflagsprop} arises.
\subsection{The space of whips}
The space of whips was analyzed in detail in the third author's work \cite{preston2011motion,preston2012geometry}. The one
modification here is the condition $\boldsymbol{\eta}_x(t,0)=\hat{\imath}$ along with $\boldsymbol{\eta}(t,0)=0$, which corresponds to holding the handle of the whip at a fixed location and orientation, rather than just at a fixed location. In this section we will derive the geodesic equation for this situation and show that the space of whips can be totally geodesically embedded into the space of regular upturned flags.

\begin{proposition}[Geodesic equation for whips]\label{whipsprop}
The geodesic equation for the Riemannian metric \eqref{whipriemannian} for curves $\boldsymbol{\eta}$ subject to $\boldsymbol{\eta}(t,0)=0$ and $\boldsymbol{\eta}_x(t,0)=\hat{\imath}$, with $\lvert \boldsymbol{\eta}_x(t,x)\rvert \equiv 1$, is given by
\begin{equation}\label{whipequation}
\frac{\partial^2\boldsymbol{\eta}}{\partial t^2} = \frac{\partial}{\partial x}\left( \sigma \frac{\partial\boldsymbol{\eta}}{\partial x}\right), \quad \boldsymbol{\eta}(t,0) = 0, \quad \boldsymbol{\eta}_x(t,0) = \hat{\imath},
\end{equation}
where $\sigma(t,x)$ is a function determined by the spatial ODE
\begin{equation}\label{whipsigma}
\frac{\partial^2\sigma}{\partial x^2} - \left\lvert \frac{\partial^2\boldsymbol{\eta}}{\partial x^2}\right\rvert^2 \sigma = -\left\lvert \frac{\partial^2\boldsymbol{\eta}}{\partial t\partial x}\right\rvert^2, \quad \sigma_x(t,0) = 0, \quad \sigma(t,1) = 0.
\end{equation}
Here $\boldsymbol{\eta}$ is assumed to be odd through $x=0$, while $\sigma$ is even through $x=0$.
\end{proposition}

\begin{proof}
The variation condition on the kinetic energy is that
\begin{equation}\label{whipvariationcondition}
\int_0^1 \langle \ddot{\boldsymbol{\eta}}(x), \varfield(x)\rangle \, dx = 0 \qquad \text{for every variation field $\varfield$.}
\end{equation}
Differentiating the equation $\langle \boldsymbol{\eta}'(x), \boldsymbol{\eta}'(x)\rangle \equiv 1$ with respect to the variation parameter, we find that every variation field must satisfy $\langle \boldsymbol{\eta}'(x), \varfield'(x)\rangle \equiv 0$. In addition we must have $\varfield(0)=0$ and $\varfield'(0)=0$. Since $\varfield'$ is orthogonal to $\boldsymbol{\eta}'$, and since $\boldsymbol{\eta}$ remains a planar curve by the whip assumption, we consider only variations
that are also planar. Hence we may write $\varfield'(x) = \delta(x) \normal(x)$ for some function $\delta$ satisfying $\delta(0)=0$, and so the variation field itself is
$$ \varfield(x) = \int_0^x \delta(s) \normal(s) \, ds, \qquad \delta(0)=0.$$
Plugging this formula into \eqref{whipvariationcondition} gives, via interchanging order of integration and then switching the variable names,
\begin{align*}
\int_0^1 \int_0^x \langle \ddot{\boldsymbol{\eta}}(x), \delta(s) \normal(s) \rangle \, ds \, dx &= \int_0^1 \delta(s) \int_s^1 \langle \ddot{\boldsymbol{\eta}}(x), \normal(s)\rangle \, ds \, dx\\
&= \int_0^1 \delta(x) \left\langle \normal(x), \int_x^1 \ddot{\boldsymbol{\eta}}(s) \, ds\right\rangle \, dx,
\end{align*}
which must be zero for every function $\delta\colon [0,1]\to\mathbb{R}$.

We find therefore that
$$\int_x^1 \ddot{\boldsymbol{\eta}}(s)\,ds = -\sigma(x) \boldsymbol{\eta}'(x)$$
for some function $\sigma$, which must satisfy $\sigma(1)=0$. Differentiating with respect to $x$ then gives
\begin{equation}\label{whipdots}
\ddot{\boldsymbol{\eta}}(x) = \frac{d}{dx} \big( \sigma(x)\boldsymbol{\eta}'(x)\big),
\end{equation}
which is equation \eqref{whipequation}. The function $\sigma$ is now a Lagrange multiplier for the condition $\lvert \boldsymbol{\eta}'(x)\rvert^2 \equiv 1$.
Differentiating that condition twice in time gives
$$\langle \ddot{\boldsymbol{\eta}}'(x), \boldsymbol{\eta}'(x)\rangle + \lvert \dot{\boldsymbol{\eta}}'(x)\rvert^2 = 0,$$
and plugging in \eqref{whipdots} gives
\begin{equation}\label{sigmaunsimplified}
\langle \sigma''(x) \boldsymbol{\eta}'(x) + 2 \sigma(x) \boldsymbol{\eta}''(x) + \sigma(x)\boldsymbol{\eta}'''(x), \boldsymbol{\eta}'(x)\rangle = -\lvert \dot{\boldsymbol{\eta}}'(x)\rvert^2.
\end{equation}
Recalling again that $\lvert \boldsymbol{\eta}'(x)\rvert^2 = 1$, successive differentiations in $x$ give $\langle \boldsymbol{\eta}'(x), \boldsymbol{\eta}''(x)\rangle = 0$
and $$\langle \boldsymbol{\eta}'(x), \boldsymbol{\eta}'''(x)\rangle + \lvert \boldsymbol{\eta}''(x)\rvert^2 = 0,$$
so that \eqref{sigmaunsimplified} becomes \eqref{whipsigma}. The boundary condition $\sigma(1)=0$ follows from the discussion above. The boundary condition at $x=0$ follows from the fact that we want $\ddot{\boldsymbol{\eta}}(0)=0$, and the compatibility condition is thus
$$\sigma'(0) \boldsymbol{\eta}'(0) + \sigma(0) \boldsymbol{\eta}''(0) = 0.$$
The inner product of this condition with $\boldsymbol{\eta}'(0)$ implies that $\sigma'(0)=0$. And while $\sigma(0)$ may not be zero, if $\boldsymbol{\eta}$ is
odd in $x$ then $\boldsymbol{\eta}''(0)=0$, which produces compatibility.

If $\boldsymbol{\eta}$ and $\dot{\boldsymbol{\eta}}$ are odd in $x$, then $\lvert \boldsymbol{\eta}''\rvert^2$ and $\lvert \dot{\boldsymbol{\eta}}'\rvert^2$ are both even in $x$,
so that $\sigma$ is even in $x$ as well. And as long as $\sigma$ remains even in $x$, equation \eqref{whipequation} ensures that $\boldsymbol{\eta}$
will remain odd in $x$. These conditions make the boundary conditions $\boldsymbol{\eta}(t,0)=0$ and $\sigma_x(t,0)=0$ redundant.
\end{proof}

At $x=0$ the compatibility condition that $\boldsymbol{\eta}(t,0)=0$, which should imply $\boldsymbol{\eta}_{tt}(t,0)=0$, requires that $\sigma_x(t,0)=0$ and that $\sigma(t,0) \boldsymbol{\eta}_{xx}(t,0)=0$. The first condition, together with $\sigma(t,1)=0$, uniquely determines the solution $\sigma$ of the ODE \eqref{whipsigma},
so the second condition cannot also be imposed. However it is satisfied automatically if $\boldsymbol{\eta}$ is assumed to be the restriction of a $C^2$ odd function.
Ensuring this compatibility is the main reason the oddness condition is convenient.

Equation \eqref{whipequation} is a nonlinear wave equation for $\boldsymbol{\eta}$, with tension determined nonlocally. Ordinarily one would specify two boundary conditions, one at $x=0$ and one at $x=1$. The fact that $\sigma(t,1)=0$ means that the natural boundary condition at $x=1$ for the symmetric differential operator $f\mapsto \tfrac{\partial}{\partial x}(\sigma \tfrac{\partial f}{\partial x})$ is simply that $f(1)$ is finite. Here the fact
that $\boldsymbol{\eta}_x(t,1)$ must be a unit vector obviates any finiteness condition, and so it is more natural to impose two conditions at $x=0$. The well-posedness theory needs to be constructed manually in any case, as no general theory applies to degenerate, nonlocal, nonlinear wave equations. See \cite{preston2011motion} and \cite{csengul2017generalized} for two approaches.

\begin{corollary}\label{thetawhipcorollary}
In terms of the function $\kappa$ defined by Lemma \ref{frenetserretlemma} and the function $\omega$ defined by
Proposition \ref{tangentspaceprop}, the equations \eqref{whipequation}--\eqref{whipsigma} take the form
\begin{alignat}{3}
\omega_t(t,x) &= \sigma(t,x) \kappa_x(t,x) + 2\sigma_x(t,x) \kappa(t,x), &\quad \omega(t,0) &= 0 \label{omegat} \\
\kappa_t(t,x) &= \omega_x(t,x), &\quad\quad \kappa(t,0) &= 0. \label{kappat} \\
\sigma_{xx}(t,x) &= \kappa(t,x)^2 \sigma(t,x) -\omega(t,x)^2, &\qquad \sigma_x(t,0) &= \sigma(t,1) = 0.\label{sigmakappaomega}
\end{alignat}
\end{corollary}

\begin{proof}
Differentiating \eqref{whipequation} with respect to $x$ gives
\begin{equation}\label{whipderiv}
\boldsymbol{\eta}_{ttx} = \sigma \boldsymbol{\eta}_{xxx} + 2\sigma_x \boldsymbol{\eta}_{xx} + \sigma_{xx} \boldsymbol{\eta}_x.
\end{equation}
Equation \eqref{etae1def} gives $\boldsymbol{\eta}_x(t,x) = \big( \cos{\theta(t,x)}, 0, \sin{\theta(t,x)}\big)$
in terms of a function $\theta$ satisfying $\theta(t,0)=0$, and if $\boldsymbol{\eta}$ is odd then additionally $\theta_x(t,0)=0$.
Plugging into equation \eqref{whipderiv} gives the components
\begin{align*}
-\theta_{tt} \sin{\theta} - \theta_t^2 \cos{\theta} &= -\sigma \theta_{xx} \sin{\theta} - \sigma \theta_x^2 \cos{\theta} - 2\sigma_x \theta_x \sin{\theta} + \sigma_{xx} \cos{\theta} \\
\theta_{tt} \cos{\theta} - \theta_t^2 \sin{\theta} &= \sigma \theta_{xx} \cos{\theta} - \sigma \theta_x^2 \sin{\theta} + 2\sigma_x \theta_x \cos{\theta} + \sigma_{xx} \sin{\theta},
\end{align*}
and resolving these gives
\begin{align}
\theta_{tt} &= \sigma \theta_{xx} + 2\sigma_x \theta_x, \label{whipthetawave} \\
-\theta_t^2 &= \sigma_{xx} -\theta_x^2 \sigma. \label{whipsigmatheta}
\end{align}

By definition of $\kappa$, we have $\boldsymbol{\eta}_{xx}(t,x) = \kappa(t,x) \normal(t,x)$ with $$\normal(t,x) = -e_2(t,x) = \big(-\sin{\theta(t,x)}, 0, \cos{\theta(t,x)}\big),$$
using Lemma \ref{frenetserretlemma} together with the fact from Remark \ref{whipremark1} that $\phi$ and $\beta$ are both zero.
Thus $\kappa(t,x) = \theta_x(t,x)$. Similarly using Remark \ref{whipremark3}, we have $\omega(t,x) = \theta_t(t,x)$. Since $\theta$ is even by Proposition \ref{oddflagsprop}, we have $\kappa(t,0)=0$ for all $t$, and since $\theta(t,0)=0$ for all $t$, we must have $\omega(t,0)=0$ for compatibility.

Thus equation \eqref{whipthetawave} implies \eqref{omegat}, while the compatibility condition $\theta_{tx} = \theta_{xt}$ implies \eqref{kappat}, and the equation \eqref{sigmakappaomega} is simply \eqref{whipsigmatheta} (which is the same as \eqref{whipsigma}) written in terms of $\kappa$ and $\omega$.
\end{proof}

We will now connect this special case back to the general case by showing that the space of whips is totally geodesic in the space of all upturned/balanced flags.
That is, we suppose that at some time we have that the flag is in the shape of a whip (i.e., $\alpha\equiv 0$) and that its velocity field will keep it that way (i.e., $\dot{\alpha}\equiv 0$). We want to prove that under these assumptions, $\ddot{\alpha}\equiv 0$ as well.

\begin{proposition}\label{whipstotallygeodesic}
The space of whips is totally geodesic in the space of upturned/balanced flags, i.e., if $\kappa(t,x)$ and $\alpha(t,x)$ solve the system in
Theorem \ref {balancedgeodesicthm}, and if $\alpha(0,x)\equiv 0$ and $\alpha_t(0,x)\equiv 0$, then $\alpha(t,x)\equiv 0$ for all $t$ and $x$.
\end{proposition}

\begin{proof}
Recall by Remark \ref{whipremark1}, the space of whips is embedded in the space of flags via the condition $\alpha\equiv 0$. Furthermore by Remark \ref{whipremark2}, the tangent space to the subspace of whips is characterized by $\dot{\alpha}\equiv 0$. Finally by Remark \ref{whipremark3}, the fact that $\alpha$ and $\dot{\alpha}$ are both zero implies that $\psi$ and $\chi$ are both identically zero. In addition we have $\gamma\equiv 1$ and $\lambda\equiv 1$, for all $x\in [0,1]$.

Under these circumstances the system \eqref{balsigma}--\eqref{xieq} simplifies to
\begin{align*}
\sigma' &= z + \frac{q}{2}, & \quad \sigma(1)&=0   \\
\rho' &= \varphi, & \rho(1) &= 0  \\
g' &= -\rho + \frac{z}{2} + \frac{q}{3}, & g(1)&= 0   \\
z' &= \kappa^2 \xi - \omega^2, & z(0) &= 0   \\
\varphi' &= q, & \varphi(0)&=0   \\
q' &= \kappa \mu, & q(0)&= 0,
\end{align*}
where $\mu = 12 g - 6\sigma$ and $\xi = 4\sigma - 6 g$. Define $G = g-\frac{1}{2} \sigma$; then the system becomes
\begin{equation}\label{reducedsystemwhip}
\sigma' = z + \frac{q}{2}, \qquad z' = \kappa^2(\sigma - 6G) - \omega^2, \quad \sigma(1)=0, \quad z(0)=0.
\end{equation}
with the other equations becoming
\begin{alignat*}{4}
\rho' &= \varphi, &\qquad G' &= -\rho + \frac{q}{12}, &\qquad \rho(1)&=0, &\quad G(1)&= 0 \\
\varphi' &= q, &\qquad q' &= 12 \kappa G, & \varphi(0)&=0, & q(0)&=0.
\end{alignat*}
These latter equations form a homogeneous system of four ODEs with solution
$$\rho\equiv G \equiv \varphi \equiv q \equiv 0.$$
The solution is unique since for any solution we have
\begin{align*}
\int_0^1 \left( 12\kappa G^2 + \frac{q^2}{12} + \varphi^2 \right)\,dx &= \int_0^1 \big( Gq' + q(G'+\rho) + \varphi \rho' \big) \,dx \\
&= \int_0^1 \frac{d}{dx} \big( Gq + \varphi \rho\big) \, dx \\
&= G(1)q(1) + \varphi(1)\rho(1) - G(0)q(0) - \varphi(0)\rho(0) \\
&= 0,
\end{align*}
using the boundary conditions. We conclude that $G$, $q$, and $\varphi$ are all zero, and thus $\rho$ must be zero as well. It follows that $g=\tfrac{1}{2}\sigma$, so that $\mu=0$ and $\xi=\sigma$ in \eqref{mueq}--\eqref{xieq}

As such equation \eqref{reducedsystemwhip} becomes
$$ \sigma'' = \kappa^2\sigma - \omega^2, \quad \sigma'(0)=0, \quad \sigma(1)=0,$$
which is precisely \eqref{sigmakappaomega}.
Then equations \eqref{kappaddot}--\eqref{alphaddot} become
\begin{align*}
\frac{\partial^2 \kappa}{\partial t^2} &= \frac{\partial^2}{\partial x^2}\big( \kappa \sigma\big)
+ \frac{\partial}{\partial x}\left(\kappa \frac{\partial \sigma}{\partial x}\right), \\
\frac{\partial^2\alpha}{\partial t^2} &= 0.
\end{align*}
The first is the spatial derivative of the equation
$$ \frac{\partial^2 \theta}{\partial t^2} =
\frac{\partial}{\partial x}\left( \sigma \, \frac{\partial \theta}{\partial x}\right) + \frac{\partial \sigma}{\partial x} \, \frac{\partial \theta}{\partial x},$$
which is precisely equation \eqref{whipthetawave}. The second shows that $\ddot{\alpha}$ will remain zero as long as $\alpha$ and $\dot{\alpha}$ are zero, as claimed.
\end{proof}

\subsection{Proof of Theorem~\ref{balancedgeodesicthm}}
{With this motivation complete, we now want to actually prove that the geodesic equation on the full space of flags is given by the equations in Theorem \ref{balancedgeodesicthm}.} We will prove it using essentially the same method as for the space of whips. Therefore we first need an expression for $\flag_{tt}(t,u,v)$. From now on we will work formally as needed, assuming that $\kappa$ and $\alpha$ have as many derivatives as any computation requires, just to derive the equations.

\begin{lemma}\label{rttlemma}
Suppose $\flag\colon [0,T]\times [0,1]^2$ is a time-dependent family of regular upturned flags given by functions $\kappa,\alpha\colon [0,T]\times [0,1]\to \mathbb{R}$ such that for each $t$, the function $\alpha$ satisfies the conditions of Definition \ref{regupflagdef}.

Then the first time derivative at $t=0$, denoted by $\dot{\flag}(u,v) = \flag_t(0,u,v)$, is given by formula \eqref{tangentflag} from Proposition \ref{tangentspaceprop}, while the second time derivative at $t=0$, denoted by $\ddot{\flag}(u,v) = \flag_{tt}(0,u,v)$, is given by
\begin{equation}\label{rttformula}
\ddot{\flag}(u,v) = \ddot{\boldsymbol{\eta}}(x) + v \big( \alpha(x) \ddot{\tangent}(x) - \ddot{\binormal}(x) \big) + \frac{v^2 \kappa(x)\dot{\alpha}(x)^2}{1+v\alpha'(x)} \,\normal(x).
\end{equation}
Here $x$ is the function of $(u,v)$ given by \eqref{xcondition}, while $\ddot{\boldsymbol{\eta}}$, $\ddot{\binormal}$, and $\ddot{\tangent}$ denote the second time derivatives at $t=0$ of $\boldsymbol{\eta}$, $\binormal$, and $\tangent$ respectively, and $\dot{\alpha}$ denotes the first time derivative as in Proposition \ref{tangentspaceprop}.
\end{lemma}

\begin{proof}
We start with the formula \eqref{rtpart2}, derived in Proposition \ref{tangentspaceprop}. Differentiate \eqref{rtpart2} again with respect to $t$, and we obtain (suppressing the independent variables on the right side):
\begin{align*}
\flag_{tt}(t,u,v) &= \boldsymbol{\eta}_{tt} + v(\alpha_t \boldsymbol{\eta}_{tx} + \alpha \boldsymbol{\eta}_{ttx} - \binormal_{tt}) + \Big( \boldsymbol{\eta}_{tx} + v\big(\alpha_x \boldsymbol{\eta}_{tx} + \alpha \boldsymbol{\eta}_{txx} - \binormal_{tx}\big) \Big) \, \frac{\partial x}{\partial t} \\
&= \boldsymbol{\eta}_{tt} + v (\alpha \boldsymbol{\eta}_{ttx} - \binormal_{tt}) + \frac{v}{1+v\alpha_x} \,\Big( (1+v\alpha_x) \alpha_t \boldsymbol{\eta}_{tx} - \alpha_t \big(\boldsymbol{\eta}_{tx} + v(\alpha_x \boldsymbol{\eta}_{tx} + \alpha \boldsymbol{\eta}_{txx} - \binormal_{tx})\big)\Big) \\
&= \boldsymbol{\eta}_{tt} + v (\alpha \boldsymbol{\eta}_{ttx} - \binormal_{tt})
+ \frac{v^2\alpha_t( \binormal_{tx} - \alpha \boldsymbol{\eta}_{txx})}{1+v\alpha_x}.
\end{align*}
Taking the time derivative of the equation $\binormal_x = \alpha \boldsymbol{\eta}_{xx}$ shows that
$$ \binormal_{tx} - \alpha \boldsymbol{\eta}_{txx} = \alpha_t \boldsymbol{\eta}_{xx} = \kappa \alpha_t \normal,$$
which gives \eqref{rttformula}.
\end{proof}

To deal with these formulas more explicitly, it is convenient to have expressions for the time derivatives of the Frenet-Serret basis. We have already derived part of this in Proposition \ref{tangentspaceprop}.

\begin{lemma}\label{frenetserretderivlemma}
Suppose $\boldsymbol{\eta}$ is a time-dependent family of curves generated by functions $\kappa$ and $\alpha$ with time derivatives $\dot{\kappa}$ and $\dot{\alpha}$, as in Proposition \ref{tangentspaceprop}. Then the first and second time derivatives of the Frenet-Serret basis vectors $\{\tangent, \normal, \binormal\}$ are given in terms of the functions $\omega, \psi, \chi$ satisfying \eqref{omegaeq}--\eqref{psieq} by
\begin{equation}\label{dotfrenet}
\begin{aligned}
\dot{\tangent} &= \omega \normal + \psi \binormal,  \\
\dot{\normal} &= -\omega \tangent + \chi \binormal, 
\\
\dot{\binormal} &= -\psi \tangent - \chi \normal,  \\
\ddot{\tangent} &= -(\omega^2+\psi^2) \tangent + (\dot{\omega}-\psi \chi) \normal + (\dot{\psi}+ \omega\chi) \binormal,  \\
\ddot{\normal} &= -(\dot{\omega} + \psi\chi) \tangent - (\omega^2 + \chi^2) \normal + (\dot{\chi}-\omega \psi)\binormal,  \\
\ddot{\binormal} &= -(\dot{\psi}-\omega\chi) \tangent - (\dot{\chi}+\psi\omega) \normal - (\psi^2+\chi^2)\binormal.
\end{aligned}
\end{equation}
Here $\dot{\omega}$, $\dot{\psi}$, and $\dot{\chi}$ represent time derivatives. In addition $\dot{\chi}$ may be determined from the other functions by the formula
\begin{equation}\label{dotchi}
\kappa \dot{\chi} = \frac{d}{dx} \dot{\psi} - \dot{\kappa}\,(\chi + \omega \alpha) - \kappa(\omega\dot{\alpha} + \alpha \dot{\omega}).
\end{equation}
\end{lemma}

\begin{proof}
The formula for the variation of the tangent is the definition of the functions $\omega$ and $\psi$ as in \eqref{omegapsidef}, and the formula for the variation of the binormal is the same as \eqref{dotbinormal1} in the proof of Proposition \ref{tangentspaceprop}.
The variation of the normal then follows from orthonormality of the Frenet-Serret basis.

To get the second variation of the tangent, we take the time derivative of the above obtained formula to get
$$ \ddot{\tangent} = \dot{\omega} \normal + \dot{\psi} \binormal + \omega \dot{\normal} + \psi \dot{\binormal},$$
and use the first variation formulas to simplify the last terms:
$$ \ddot{\tangent} = \dot{\omega} \normal + \dot{\psi} \binormal
+ \omega (-\omega \tangent + \chi \binormal) + \psi (-\psi \tangent - \chi\normal)
= -(\omega^2+\psi^2) \tangent + (\dot{\omega} -\psi\chi)\normal + (\dot{\psi} + \omega\chi) \binormal.$$
The remaining second variation formulas are obtained the same way. Finally, equation \eqref{dotchi} comes from differentiating the formula \eqref{chieq} with respect to $t$, and solving for $\dot{\chi}$.
\end{proof}

%
%
%
%
\begin{proof}[Proof of Theorem~\ref{balancedgeodesicthm}]
The condition to be a geodesic is that the acceleration is perpendicular to all variation vectors, as in \eqref{variationalacceleration}. Upon changing to $(x,y)$ coordinates and using the change-of-variables formula \eqref{jacobianarea}, this condition becomes
\begin{equation}\label{variationalaccelerationxy}
\int_0^1 \int_0^{\gamma(x)} \big(1+y\alpha'(x)\big)\langle \ddot{r}(x,y), \varfield(x,y)\rangle \, dy\,dx = 0, \qquad \text{for every variation field $\varfield$.}
\end{equation}
Using Proposition \ref{tangentspaceprop}, we may write a general perturbation vector as
\begin{equation}\label{tangentperturbation}
\varfield = \int_0^x \big( \delta(s) \normal(s) + \varepsilon(s) \binormal(s)\big)\,ds + y \Big(
\frac{\varepsilon'(x)}{\kappa(x)} \normal(x) + \varepsilon(x)\big(\tangent(x) + \alpha(x)\binormal(x)\big)\Big),
\end{equation}
for some functions $\delta$ and $\varepsilon$ satisfying $\delta(0)=\varepsilon(0)=\varepsilon'(0)=0$, but which are otherwise arbitrary.
Using Lemma \ref{rttlemma}, we have in the $(x,y)$ coordinates that
\begin{multline}\label{rttformulajacob}
\big(1+y\alpha'(x)\big) \ddot{\flag}(x,y) = \big(1+y\alpha'(x)\big) \big(\ddot{\boldsymbol{\eta}}(x) + y \nu(x)\big) + y^2 \pi(x)\\
\text{where } \nu(x) =  \alpha(x) \ddot{\tangent}(x) - \ddot{\binormal}(x), \quad \text{and}\quad
\pi(x) = \kappa(x)\dot{\alpha}(x)^2 \,\normal(x).
\end{multline}
Hence the condition, after performing the integration of the powers of $y$ for $0\le y\le \gamma(x)$, is that
\begin{multline}\label{integralyzeroacc}
\int_0^1 \int_0^x \Big\langle \delta(s) \normal(s) + \varepsilon(s) \binormal(s)\big), \\
\gamma(x)\ddot{\boldsymbol{\eta}}(x) + \tfrac{\gamma(x)^2}{2} \big(\nu(x)+\alpha'(x)\ddot{\boldsymbol{\eta}}(x)\big) + \tfrac{\gamma(x)^3}{3} \big( \pi(x) + \alpha'(x) \nu(x)\big)   \Big\rangle \, ds\, dx \\
+ \int_0^1 \Big\langle \tfrac{1}{\kappa(x)} \varepsilon'(x) \normal(x) + \varepsilon(x)\big( \tangent(x) + \alpha(x)\binormal(x)\big), \\
\tfrac{\gamma(x)^2}{2} \ddot{\boldsymbol{\eta}}(x) + \tfrac{\gamma(x)^3}{3} \big( \nu(x) + \alpha'(x)\ddot{\boldsymbol{\eta}}(x)\big) + \tfrac{\gamma(x)^4}{4} \big( \pi(x) + \alpha'(x)\nu(x)\big) \Big\rangle \, dx = 0
\end{multline}

Notice that $\alpha'(x)$ appears here only in combination with $\gamma(x)$; thus recalling the formula $\lambda(x) = 1+\alpha'(x)\gamma(x)$ from \eqref{lambdagammadef}, we can define new functions $K, J\colon [0,1]\to \mathbb{R}$ to be the second terms in the formula above, by the formulas
\begin{align}
K(x) &= -\int_x^1 \Big( \frac{\gamma(s)(1+\lambda(s))}{2} \ddot{\boldsymbol{\eta}}(s) + \frac{\gamma(s)^2(1+2\lambda(s))}{6} \nu(s) + \frac{\gamma(s)^3}{3} \pi(s)\Big) \, ds \label{Kdef}\\
J(x) &= \frac{\gamma(x)^2(1+2\lambda(x))}{6} \ddot{\boldsymbol{\eta}}(x) + \frac{\gamma(x)^3(1+3\lambda(x))}{12} \nu(x) + \frac{\gamma(x)^4}{4} \pi(x).\label{Jdef}
\end{align}
Then we obviously have $K(1)=0$, along with
\begin{equation}\label{Kprime}
K'(x) = \frac{\gamma(x)(1+\lambda(x))}{2} \ddot{\boldsymbol{\eta}}(x) + \frac{\gamma(x)^2(1+2\lambda(x))}{6} \nu(x) + \frac{\gamma(x)^3}{3} \pi(x).
\end{equation}

In terms of $J$ and $K$, formula \eqref{integralyzeroacc} becomes
\begin{multline}
-\int_0^1 \langle \delta(x) \normal(x) + \varepsilon(x) \binormal(x)\big),
K(x)\rangle + \Big\langle \frac{1}{\kappa(x)} \varepsilon'(x) \normal(x) + \varepsilon(x)\big( \tangent(x) + \alpha(x)\binormal(x)\big), J(x)
\Big\rangle \, dx = 0.
\end{multline}
Now integrate by parts to remove the derivative on $\varepsilon$, and we obtain
\begin{multline}\label{finalvariationconditionbalanced}
-\int_0^1 \delta(x) \langle \normal(x), K(x)\rangle \, dx
+ \frac{\varepsilon(x)}{\kappa(x)} \langle J(x), \normal(x)\rangle\Big|_{x=0}^{x=1} \\
+ \int_0^1 \varepsilon(x) \left[ -\langle \binormal(x), K(x)\rangle - \frac{d}{dx} \left( \frac{\langle J(x), \normal(x)\rangle}{\kappa(x)}\right) + \langle J(x), \tangent(x)+\alpha(x)\binormal(x) \rangle \right]\,dx = 0.
\end{multline}
This expression must be zero for every choice of functions $\delta$ and $\varepsilon$. We conclude that
\begin{align}
\langle K(x), \normal(x)\rangle &\equiv 0 \label{noKnormal} \\
\Big\langle \frac{J(1)}{\kappa(1)}, \normal(1)\Big\rangle &= 0, \qquad K(1)=0, \label{Jendpoint} \\
\frac{d}{dx}\left( \frac{\langle J(x), \normal(x)\rangle}{\kappa(x)}\right) &= - \langle \binormal(x), K(x)\rangle + \langle J(x), \tangent(x)+\alpha(x)\binormal(x)\rangle.  \label{gequationJ}
\end{align}

Now consider the five vectors here:
$J$, $K$, $\ddot{\boldsymbol{\eta}}$, $\ddot{\tangent}$, and $\nu$. By Lemma \ref{frenetserretderivlemma}, the terms $\ddot{\tangent}$, $\ddot{\normal}$, and $\ddot{\binormal}$ (and thus $\nu$ from \eqref{rttformulajacob}) all depend on only three functions $\{\dot{\omega}, \dot{\psi}, \dot{\chi}\}$, which are unknown, and finding them will yield the geodesic equation.
Also $K$ has two unknown components by \eqref{noKnormal}, while $J$ and $\ddot{\boldsymbol{\eta}}$ have three: this is eleven unknowns in all.
Formula \eqref{dotchi} gives a differential equation for $\dot{\chi}$ in terms of $\dot{\psi}$, so we have ten unknowns left. Then the fact that $\frac{d}{dx} \ddot{\boldsymbol{\eta}} = \ddot{\tangent}$ gives a differential equation that solves for three of these components, leaving us with seven unknowns. Equation \eqref{gequationJ} gives another differential equation, leaving us with six unknowns. The components of equations \eqref{Jdef} and \eqref{Kprime} form another six equations (two of which are differential equations), which should in principle completely solve for all unknowns. So we have six linear differential equations and five linear algebraic equations for the eleven unknowns. Three boundary conditions for the six differential equations at $x=1$ are given by the three conditions \eqref{Jendpoint}, while three more are given at $x=0$ by the requirements that $\ddot{\boldsymbol{\eta}}$ and $\nu$ are all zero (and the fact that $\kappa(0)=0$ means we need no more than three conditions here).

We now write $J$, $K$, $\ddot{\boldsymbol{\eta}}$, and $\nu$ in the Frenet-Serret basis, defining eleven functions $f$, $g$, $h$, $\sigma$, $\rho$, $z$, $\xi$, $\varphi$, $p$, $q$, and $\mu$ as
\begin{align}
J(x) &= f(x)\tangent + \kappa(x) g(x) \normal(x) + h(x) \binormal(x), \label{Jcomponents}\\
K(x) &= \sigma(x)\tangent(x) + \rho(x) \binormal(x), \label{Kcomponents} \\
\ddot{\boldsymbol{\eta}}(x) &= z(x) \tangent(x) + \kappa(x) \xi(x) \normal(x) + \varphi(x) \binormal(x), \label{ddotetacomponents}\\
\ddot{t}(x) &= -\big( \omega(x)^2 + \psi(x)^2\big) \tangent + p(x) \normal(x) + q(x) \binormal(x), \label{tddotcomponents}\\
\kappa(x) \mu(x)&= \langle \nu(x), \normal(x)\rangle. \label{mudef}
\end{align}
where $K$ has no normal component by equation \eqref{noKnormal}.
Here $\mu$, $\xi$, $f$, $h$, and $p$ will be determined algebraically in terms of the fundamental variables $(g,\sigma,\rho,z,\varphi,q)$, the latter of which solve six coupled one-dimensional ODEs.

By the second variation formula for the tangent vector in Lemma \ref{frenetserretderivlemma}, the equation for $\ddot{\boldsymbol{\eta}}(x)$ becomes
\begin{equation}\label{pqdots}
p(x) = \dot{\omega}(x) - \psi(x)\chi(x), \qquad \text{and}\qquad
q(x) = \dot{\psi}(x) + \omega(x) \chi(x).
\end{equation}
Next we use Lemma \ref{frenetserretderivlemma} again to get
\begin{equation}\label{nuexplicit}
\begin{split}
\nu &= \alpha \ddot{\tangent} - \ddot{\binormal}\\
&= \big( \dot{\psi}-\omega\chi-\alpha(\omega^2+\psi^2)\big)\tangent +
     \big(  \alpha (\dot{\omega}-\psi \chi) + \dot{\chi} + \psi \omega\big) \normal
     + \big( \alpha(\dot{\psi}+\omega\chi) + \psi^2+\chi^2 \big) \binormal.
\end{split}
\end{equation}
We use equation \eqref{dotchi} to rewrite the term $\dot{\chi}$ appearing above in terms of $p$ and $q$ defined by \eqref{pqdots}, to obtain
\begin{equation}\label{dotchi2}
\begin{split}
\kappa \dot{\chi} &= \frac{d}{dx} (q - \omega\chi)
- \dot{\kappa}(\chi + \omega \alpha) - \kappa\omega\dot{\alpha}
-\kappa \alpha (p+\psi\chi) \\
&= q' + \omega(\kappa\dot{\alpha} + \alpha \dot{\kappa} + \kappa\psi)
- \chi (\dot{\kappa} - \kappa\alpha\psi)
- \dot{\kappa}(\chi + \omega \alpha) - \kappa\omega\dot{\alpha}
-\kappa \alpha (p+\psi\chi)  \\
&= q' - \kappa \alpha p
+ \omega\kappa\psi - 2 \chi \dot{\kappa},
\end{split}
\end{equation}
after using \eqref{omegaeq}  and \eqref{psieq}.

Using this, we find that the normal component of $\nu$ from \eqref{nuexplicit} is
\begin{equation}\label{normalnu}
\kappa^2 \mu = \kappa \langle \nu, \normal\rangle = q' - 2 \chi \dot{\kappa}  + 2 \omega\kappa \psi,
\end{equation}
which is equation \eqref{balq}. The condition $q(0)=0$ comes from \eqref{pqdots} and the requirement that $\dot{\psi}(0)=0$ and $\chi(0)=0$,
from equations \eqref{chieq}--\eqref{psieq}.

To find $\mu$ and $\xi$, we solve the system coming from the normal components of \eqref{Jdef}--\eqref{Kprime}, which becomes using \eqref{Jcomponents}--\eqref{ddotetacomponents} the system
\begin{equation}\label{muxisystem}
\begin{split}
g &= \frac{\gamma^2(1+2\lambda)}{6} \,\xi+ \frac{\gamma^3(1+3\lambda)}{12} \,\mu + \frac{\gamma^4}{4}\, \dot{\alpha}^2,\\
\sigma+\alpha \rho &=  \frac{\gamma(1+\lambda)}{2} \,\xi + \frac{\gamma^2(1+2\lambda)}{6} \,\mu + \frac{\gamma^3}{3} \,\dot{\alpha}^2,
\end{split}
\end{equation}
which can be solved for $\mu$ and $\xi$ to obtain \eqref{mueq}--\eqref{xieq}.

Using the components \eqref{ddotetacomponents}--\eqref{tddotcomponents}, the condition $\frac{d}{dx} \ddot{\boldsymbol{\eta}} = \ddot{\tangent}$ yields the equations
$$ (z' - \kappa^2 \xi) \tangent + \Big(\frac{d}{dx}(\kappa \xi) + \kappa(z+\alpha \varphi)\Big)\normal
+ (\varphi' - \kappa^2 \alpha \xi) \binormal
= -(\omega^2 + \psi^2) \tangent + p \normal + q \binormal
$$ 
The tangent and binormal components of this are equations \eqref{balz}--\eqref{balphi}, and the conditions that $z(0)=\varphi(0)=0$ come from the fact that $\ddot{\boldsymbol{\eta}}(0)=0$.
The normal component is
\begin{equation}\label{peqn}
p = \frac{d}{dx}(\kappa \xi)  + \kappa(z+\alpha\varphi),
\end{equation}
which we will return to in a moment.

Conditions \eqref{Jendpoint} and \eqref{gequationJ} imply that
\begin{equation}\label{geqfhrho}
g'(x) = -\rho(x) + f(x)+\alpha(x)h(x), \qquad g(1)=0.
\end{equation}
Using the definition \eqref{Jdef} of $J$ and the formulas \eqref{Jcomponents} and \eqref{ddotetacomponents}, as well as the explicit formulas for $\nu$ from \eqref{nuexplicit} and the definition \eqref{pqdots} of $q$, we get for the tangential and binormal components that
\begin{align*}
f &= \frac{\gamma^2(1+2\lambda)}{6} \, z + \frac{\gamma^3(1+3\lambda)}{12} \, \big(q - 2\omega\chi - \alpha(\omega^2 + \psi^2)\big), \\
h &= \frac{\gamma^2(1+2\lambda)}{6} \, \varphi + \frac{\gamma^3(1+3\lambda)}{12} \, ( \alpha q + \psi^2+\chi^2),
\end{align*}
and substituting these formulas into \eqref{geqfhrho} gives \eqref{balg}.

Finally we use the formula \eqref{Kprime} for $K'$ along with the component definition \eqref{Kcomponents}, together with \eqref{ddotetacomponents} and \eqref{nuexplicit}, to get
\begin{align*}
\sigma' &= \frac{\gamma(1+\lambda)}{2} \, z + \frac{\gamma^2(1+2\lambda)}{6} \,
\big( q - 2\omega\chi - \alpha(\omega^2 + \psi^2)\big), \\
\rho' &= \frac{\gamma(1+\lambda)}{2} \, \varphi + \frac{\gamma^2(1+2\lambda)}{6} \, \big(\alpha q + \psi^2+\chi^2 \big),
\end{align*}
which are exactly \eqref{balsigma}--\eqref{balrho}. The boundary conditions $\sigma(1)=\rho(1)=0$ come from $K(1)=0$, by construction \eqref{Kdef} of $K$.

We now have a closed system \eqref{balsigma}--\eqref{balq} of six first-order  differential equations with six boundary conditions, and its solution will give us $p$ and $q$; we get $q$ directly, while $p$ comes from \eqref{peqn} since we know $\xi$. These will give us $\dot{\omega}$ and $\dot{\psi}$ by formulas \eqref{pqdots}. It remains to see how these tell us the second derivatives of the basic functions $\kappa$ and $\alpha$. To get these, we differentiate the equations \eqref{omegaeq}--\eqref{psieq} with respect to $t$, and solve for $\ddot{\kappa}$ and $\ddot{\alpha}$.

Differentiating \eqref{omegaeq} in time, then using \eqref{pqdots}, \eqref{chieq}, and \eqref{psieq} gives
\begin{equation}\label{kappaddotpartial}
\begin{split}
\ddot{\kappa} &= \dot{\omega}' + \dot{\kappa} \alpha \psi + \kappa \psi \dot{\alpha} + \kappa \alpha \dot{\psi} \\
&= \frac{d}{dx}(p+\psi \chi) + \dot{\kappa} \alpha \psi + \kappa \psi \dot{\alpha} + \kappa \alpha (q - \omega \chi) \\
&= p' - \psi(\kappa \dot{\alpha}+\alpha\dot{\kappa}+\kappa\psi) +
\chi (\kappa \chi + \kappa \alpha \omega)
+ \dot{\kappa} \alpha \psi + \kappa \psi \dot{\alpha} + \kappa \alpha (q - \omega \chi) \\
&= p' + \kappa \alpha q + \kappa(\chi^2 - \psi^2).
\end{split}
\end{equation}
Finally recalling \eqref{peqn} and plugging it in for $p'$, we get \eqref{kappaddot} for $\ddot{\kappa}$.

To get $\ddot{\alpha}$, we do the same thing: differentiate \eqref{psieq} in time to get
\begin{equation}\label{alphaddotraw}
0 = \kappa\ddot{\alpha} + 2\dot{\kappa}\dot{\alpha} + \alpha\ddot{\kappa}
+ \psi \dot{\kappa} + \kappa \dot{\psi} + \frac{d}{dx} \dot{\chi}.
\end{equation}
We use equation \eqref{kappaddotpartial} to eliminate $\ddot{\kappa}$, and note that equations \eqref{dotchi2} and \eqref{normalnu} now imply
$$ \dot{\chi} = \kappa \mu - \omega \psi - \alpha p.$$
Hence equation \eqref{alphaddotraw} becomes, using \eqref{omegaeq}--\eqref{chieq},
\begin{align*}
0 &= \kappa \ddot{\alpha} + 2\dot{\kappa}\dot{\alpha} +
\alpha \big( p' + \kappa \alpha q + \kappa(\chi^2 - \psi^2)\big)
+ \psi\dot{\kappa} + \kappa \dot{\psi} +
\frac{d}{dx}\big( \kappa \mu - \omega \psi - \alpha p\big)  \\
&= \kappa\ddot{\alpha} + 2\dot{\kappa}\dot{\alpha}
- \alpha' p + \kappa \alpha^2 q + \kappa\alpha (\chi^2 - \psi^2)
+ \psi \dot{\kappa} + \kappa (q-\omega\chi)
+ \frac{d}{dx}(\kappa \mu)  - \omega\psi' - \psi \omega' \\
&= \kappa\ddot{\alpha} + 2\dot{\kappa}\dot{\alpha}
- \alpha' p + \kappa (1+\alpha^2) q + \kappa\alpha (\chi^2 - \psi^2)
+ \psi \dot{\kappa} - \kappa\omega\chi
+ \frac{d}{dx}(\kappa \mu)  \\
&\qquad\qquad - \omega(\kappa\chi + \kappa\alpha\omega)
- \psi (\dot{\kappa}-\kappa\alpha\psi) \\
&=  \kappa\ddot{\alpha} + 2\dot{\kappa}\dot{\alpha}
- \alpha' p + \kappa (1+\alpha^2) q
+ \kappa\alpha (\chi^2 - \omega^2) - 2 \kappa \omega \chi + \frac{d}{dx}(\kappa \mu).
\end{align*}
Using \eqref{peqn} to eliminate $p$, we obtain \eqref{alphaddot}.
\end{proof}

\subsection{{Towards wellposedness of the geodesic equations}}
 {Our primary goal in this work has been to derive the geodesic equations; unfortunately we have not been able to prove local existence for them yet, and have to leave this as a goal for future research. As compared to the situation of either ideal fluids~\cite{Arn2014,EM1970} or the situation of whips (inextensible strings)~\cite{preston2011motion,preston2012geometry}  there are several additional issues here, which we will explain next. }

The first one is existence and uniqueness of the six-dimensional ODE system \eqref{balsigma}--\eqref{xieq}, needed to find the tension coefficients. In the whip case as in Proposition \ref{whipstotallygeodesic}, we could show that four of the functions are zero and compute a conservation law to obtain uniqueness of this solution, with the remaining two reducing to the tension $\sigma$ and its spatial derivative, solving the second-order tension ODE \eqref{whipsigma}. Existence and uniqueness for a second-order boundary value problem is straightforward, but for this system one would still need to establish this. Beyond that, one then has to prove local well-posedness for the nonlinear wave equations \eqref{kappaddot} and \eqref{alphaddot} for $\kappa$ and $\alpha$, where unfortunately the coefficients are not smooth (since $\lambda$ has a kink at $x^*$ in the upturned case). This is probably quite difficult. Possibly the best one can hope for is weak solutions along the lines of those found for whips by {\c{S}}eng{\"u}l and Vorotnikov~\cite{csengul2017generalized}.

{
One might believe that these difficulties have their source in our coordinate representation of the space of flags, and that the geodesic equations would be better behaved 
in their  ``natural coordinates,'' i.e., when working directly with the space of functions  $\flag:[0,1]\times[0,1]\to \mathbb R^3$, satisfying the isometry and flag pole conditions.
Using this representation and Theorem~\ref{thm:submanifold} one can view the space of (smooth) flags as a submanifold of a flat space; note that this is the analogue to Arnold's description of the Euler-equation by viewing the space of volume preserving diffeomorphisms as a submanifold of the group of all diffeomorphisms equipped with the flat $L^2$-metric~\cite{Arn2014}. Using this representation the geodesic equation can be obtained  from the general principle that geodesics in a submanifold of a flat space satisfy the condition that the acceleration is normal to the submanifold. }

{Deriving the geodesic equation from this principle relies on computing the orthogonal complement of the tangent space, which turns out to be  rather involved; mainly due to the boundary conditions at the top and bottom edges. The exact derivation of these formulas is extremely tedious and thus we refrain from presenting them here; we refer the interested reader to a previous version of the present article which is available at \url{https://arxiv.org/abs/1905.06378v1}, where we presented the derivation in detail. In the following  
we will only sketch the resulting formulas and describe the difficulties with this approach.}

{Inspired by the formula $\partial_u(A \mathbf{r}_u)$ for the orthogonal component of a vector in the case of inextensible threads, c.f. \cite{preston2011motion, preston2012geometry}, it is natural to consider the ansatz {that an orthogonal vector $\mathbf{w}$ takes the form}
\begin{equation}\label{orthocomplement}
\mathbf{w} = \partial_u(A\mathbf{r}_u + B\mathbf{r}_v) + \partial_v(B\mathbf{r}_u + C\mathbf{r}_v).
\end{equation}
It turns out that it is easy to see that any such field is indeed in the orthogonal complement, i.e., it is orthogonal to any tangent vector after imposing the correct boundary conditions} {on $A$, $B$, and $C$.} {We will omit the exact form of the boundary conditions to keep the presentation simple. 
Using this, one then obtains the geodesic equation in the form
\begin{equation}\label{flaggeodesicwave}
\mathbf{r}_{tt} = \partial_u(A \mathbf{r}_u +B\mathbf{r}_v) + \partial_v(B \mathbf{r}_u+C\mathbf{r}_v) 
\end{equation}
for three functions $A$, $B$, $C$, which act as the ``tensions'' in a nonlinear wave equation. Here the right side is essentially the second fundamental form of the space of isometric immersions inside the space of all immersions. The functions $A$, $B$, $C$ are given in terms of the velocity $\mathbf{r}_t$ in a  similar way as in the whip equation \eqref{whipequation}: they satisfy the equations
\begin{align*}
A_{uu} + B_{uv} - eK &= -\mathbf{r}_{tu}\cdot \mathbf{r}_{tu} \\
A_{uv} + B_{uu} + B_{vv} + C_{uv} - 2fK &= -\mathbf{r}_{tu}\cdot \mathbf{r}_{tv} \\
B_{uv} + C_{vv} - gK &= -\mathbf{r}_{tv}\cdot \mathbf{r}_{tv} \\
K &= eA + 2fB + gC,
\end{align*}
together with some rather complicated boundary conditions, which we again omit. The difficulty in solving this system is that the equations are again rather degenerate. However, because of this one is able to get quite far with explicit but very involved formulas for the solutions; ultimately two rather complicated differential equations for a single function determine $A$, $B$, and $C$, and thus the geodesic equation. }

{Apart from these technical difficulties, which stem mostly from the complicated boundary terms, there arise more fundamental problems for using these coordinates to obtain local well-posedness:  in the context of Arnold's description for incompressible fluids, it has been shown by Ebin and Marsden that the orthogonal projection is well behaved~\cite{EM1970}; even for functions which only admit finite differentiability. This allowed them to interpret the Euler equations as an ODE on a Banach space of functions, and consequently they were able to deduce the local wellposedness using  Picard iteration. Investigating the proof of our submanifold result in Section~\ref{submanifoldsection}, one can see that the orthogonal projection for the space of isometric immersions is not as well behaved: indeed there is always a loss of regularity  in the orthogonal projection, and thus
one cannot hope to interpret the geodesic equation as an ODE on a Banach space. Thus similar techniques as used by Ebin and Marsden for the Euler equation are  bound to fail in our situation.}

\section{Numerical Experiments}\label{sec:numerics}
In this final section we will present numerical experiments demonstrating the behavior of the developed theory. We want to emphasize that our goal here is not to develop a comprehensive numerical framework, but only to show some numerical experiments detailing the developed theory. In future work we plan to develop a more serious numerical framework for both the geodesic initial and boundary value problem. For the purpose of the current article we focus on the geodesic boundary value problem, i.e., given two upturned flags $\flag_0$ and $\flag_1$, we aim to minimize the Riemannian energy
\begin{equation}
E(\flag)=\int_0^1
\llangle \dot{\flag}(t), \dot{\flag}(t)\rrangle_{\flag(t)} dt = \int_0^1 \int_0^1 \int_0^1 \lvert \dot{\flag}(t,u,v)\rvert^2 \, du\,dv\, dt
\end{equation}
over all paths of upturned flags $\flag$ such that $\flag(0)=\flag_0$ and $\flag(1)=\flag_1$. Here we view $\flag$ both as a map from $[0,1]$ into the space as flags, as well as a function $[0,1]\times[0,1]\times[0,1]\to \mathbb R^3$. We then use the coordinates from Theorem~\ref{flagsolutionthm} to represent a flag via its functions of one variable $\alpha$ and $\kappa$. Recall that the function $\alpha$ satisfies a series of constraints \eqref{alphaconditionsagain}, which can be written as
\begin{equation}\label{alphaconstraints2}
\alpha(0)=0, \qquad \alpha(1)>0, \qquad  \max(1-x, \alpha(t,x)) +(1-x)\alpha_x(t,x)>0
\end{equation}
Consequently  paths of flags corresponds to paths of pairs of functions $\alpha(t,x)$ and $\kappa(t,x)$, where $\alpha$ satisfies for each $t$ the constraints~\eqref{alphaconstraints2}.
We then express the Riemannian energy in terms of $\alpha(t,x)$ and $\kappa(t,x)$ using Proposition~\ref{metric_coordinates}, which reduces the geodesic boundary value problem to a constrained minimization problem for  $\alpha(t,x)$ and $\kappa(t,x)$.

Next we discretize the functions $\alpha$ and $\kappa$: for fixed $t$ we  represent the functions $\alpha$ and $\kappa$ as piecewise linear functions with $N_x$ break points, i.e., we reduce $\alpha$ and $\kappa$  to their values $\alpha(t,x_i)$ and $\kappa(t,x_i)$, where $x_i=\frac{i}{N_x}$ for $i\in 0\ldots N_x$. Thereby we have discretized a flag as a (constrained) vector in $\mathbb R^{2N_x+2}$.
Similarly we discretize in time using $N_t$ time steps, which allows us to represent a path of  flags $\flag(t,u,v)$ as a (constrained) vector in $\mathbb R^{(2N_x+2)(N_t+1)}$.
Note  that $\alpha(0,x_i)$, $\kappa(0,x_i)$, $\alpha(1,x_i)$ and $\kappa(1,x_i)$ are prescribed, i.e., they are the given boundary conditions. In addition $\alpha(t_i,0)=0$. Thus in practice we have $\mathbb R^{(2N_x+1)(N_t-1)}$ free variables in our constrained minimization problem; the remaining constraints stem from the conditions on $\alpha(t_i,\cdot)$. To deal with the remaining constraints we tried two different approaches: first, we simply optimized the unconstrained minimization problem and checked subsequently if the solution satisfies the desired constraints. This strategy works well for flags that are not too far apart. If one considers, however, situations where the boundary points (flags) are too far apart, we observed that the minimizer tends to deform the path in a direction such that the constraints are indeed violated. To overcome this difficulty we
tried to relax  the remaining constraints and add them instead as a penalty function to the energy functional.  This would again allow us to tackle the an unconstrained minimization problem. This naive approach did not seem to lead to satisfactory results, and we believe that more sophisticated methods will be necessary to handle such situations. We will leave this part open for future research.

 In the following we will detail the necessary steps to implement the kinetic energy of a path of flags. We will write $A\in \mathbb R^{N_1\times \dots\times N_k}$ to denote an array of dimensions $[N_1,N_2,\ldots, N_k]$.
\begin{itemize}
\item {\bf Input:} A PL path of flags, described via its generating functions $\alpha(t_i,x_j)$ and  $\kappa(t_i,x_j)$ with $i\in 0,\ldots N_t$ and $j\in 0,\ldots N_x$; i.e.,  both $\alpha$ and $\kappa$ are arrays of dimensions $[N_t+1,N_x+1]$, where the first dimension refers to the discrete time points, and the second dimension to the discretization in space. \item For each time $t_i$, $i\in [0,N_t]$ set up $\gamma\in \mathbb R^{(N_t+1)\times (N_x+1)}$ and $\lambda\in  \mathbb R^{(N_t+1)\times (N_x+1)}$ using equation~\eqref{lambdagammadef}.
\item For each time $t_i$, $i\in [0,N_t]$ use finite differences in the second dimension of $\alpha$ to calculate the $x$-derivative $\alpha_x  \in \mathbb R^{(N_t+1)\times (N_x+1)}$.
\item For each time $t_i$, $i\in [0,N_t]$ use the Euler method to solve the ODE system~\eqref{thetaeq}-\eqref{betaeq} to calculate $(\theta,\phi,\beta)\in (\mathbb R^{(N_t+1)\times (N_x+1)},\mathbb R^{(N_t+1)\times (N_x+1)}, \mathbb R^{(N_t+1)\times (N_x+1)})$.
\item For each time $t_i$, $i\in [0,N_t]$ use equations~\eqref{etae1def} and \eqref{e2e3def} to set up the orthonormal basis
 $(\boldsymbol{\eta}_x,e_2,e_3)\in (\mathbb R^{N_t\times N_x\times 3},\mathbb R^{(N_t+1)\times (N_x+1)\times 3},\mathbb R^{(N_t+1)\times (N_x+1)\times 3})$.
 \item For each  $t_i$, $i\in [0,N_t]$ use equation ~\eqref{normalbinormaldef} to set up the binormal vector $\binormal\in  \mathbb R^{(N_t+1)\times (N_x+1)\times 3}$.
 \item For each $x_j$, $j\in [0,N_x]$ use finite differences of $\binormal$ in the first dimension to calculate $\binormal_t$.
\item For each $x_j$, $j\in [0,N_x]$ use finite differences of $\boldsymbol{\eta}_x$ in the first dimension to calculate $\boldsymbol{\eta}_{tx}$.
\item For each  $t_i$, $i\in [0,N_t]$ use integration in $x$ (the second dimension) to calculate $\boldsymbol{\eta}_t$.
\item For each  $t_i$, $i\in [0,N_t]$ calculate $\langle \dot \flag,\dot \flag\rangle_{\flag}\in \mathbb R^{N_t+1}$ using equation~\eqref{eq:metric_alt}.
\item Calculate the total Energy by summing up the array $1/N_t*\langle \dot \flag,\dot \flag\rangle_{\flag}$.
\item {\bf Return:} Energy.
 \end{itemize}
Note that for a given minimizer we can calculate the constraint minimization \eqref{alphaconstraints2} via
$$\operatorname{Cn} =\operatorname{min}\left(0,\sum_{i=1}^{N_t}\sum_{j=1}^{N_x} \operatorname{max}\{1-x_j, \alpha(t_i,x_j)\} +(1-x_j)\alpha_x(t_i,x_j)\right).$$

We implemented the energy functional in {\tt pytorch}, which allows us to take advantage of the automatic differentiation capabilities to calculate the gradient and use the L-BFGS algorithm, as introduced in \cite{liu1989limited}, to minimize the energy. To initialize the L-BFGS algorithm one needs to specify an initial guess for the solution, for which we can choose e.g., the linear interpolation between the generating functions of the given source and target flag.

\begin{figure}
\includegraphics[width=.24\textwidth]{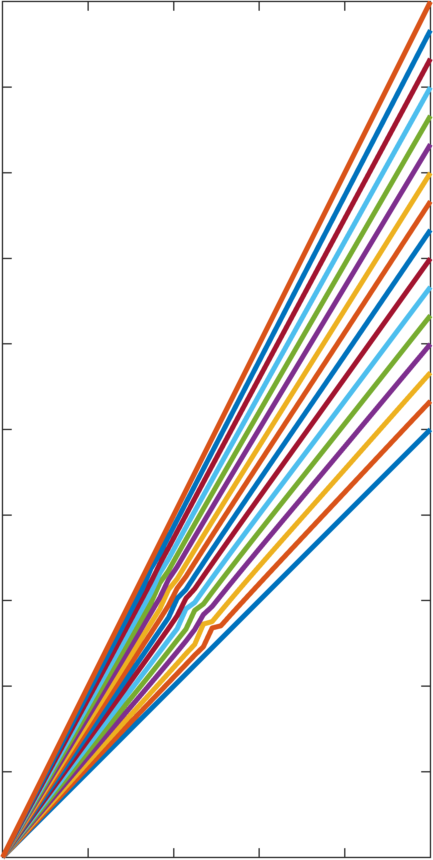}
\qquad
\includegraphics[width=.6\textwidth]{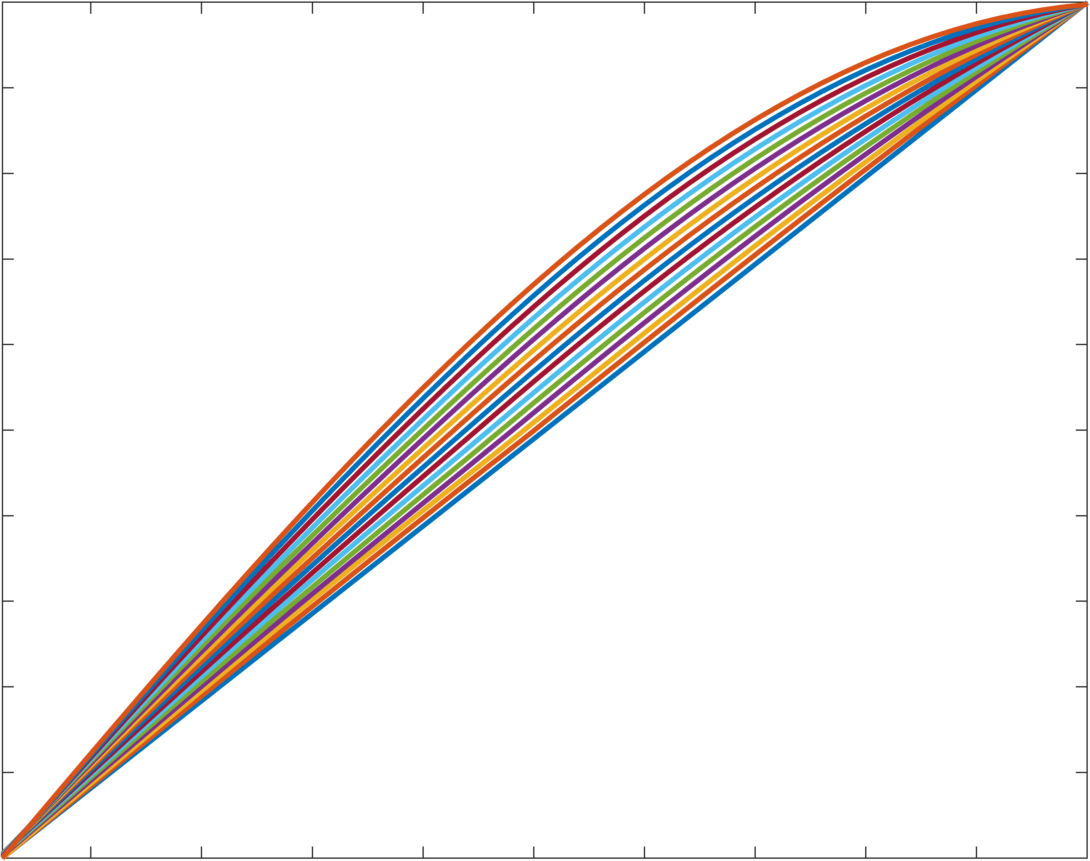}\\
\caption{A geodesic in $\alpha$-$\kappa$  space. Left figure: time evolution of $\alpha$; right figure: time evolution of $\kappa$.}
\label{fig:geod:alphaspace}
\end{figure}

\begin{figure}
\includegraphics[width=.1\textwidth]{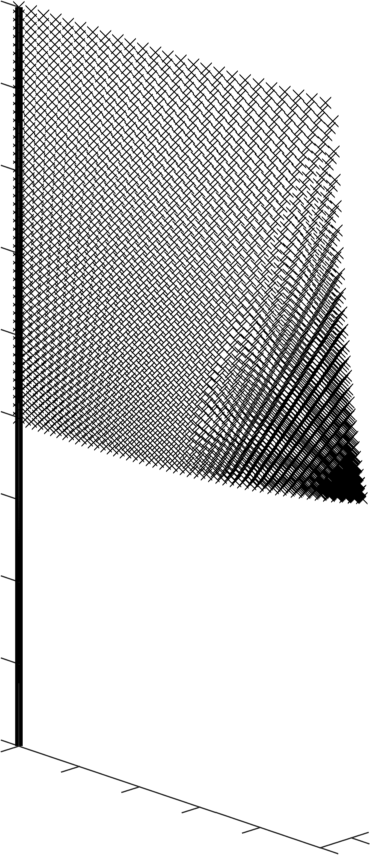}
\qquad
\includegraphics[width=.1\textwidth]{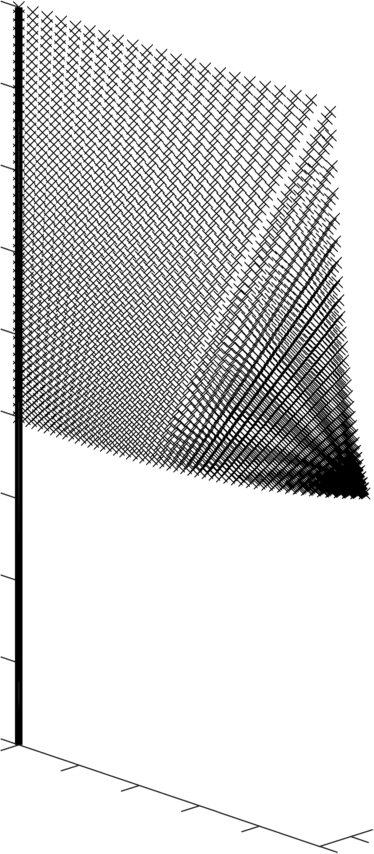}
\qquad
\includegraphics[width=.1\textwidth]{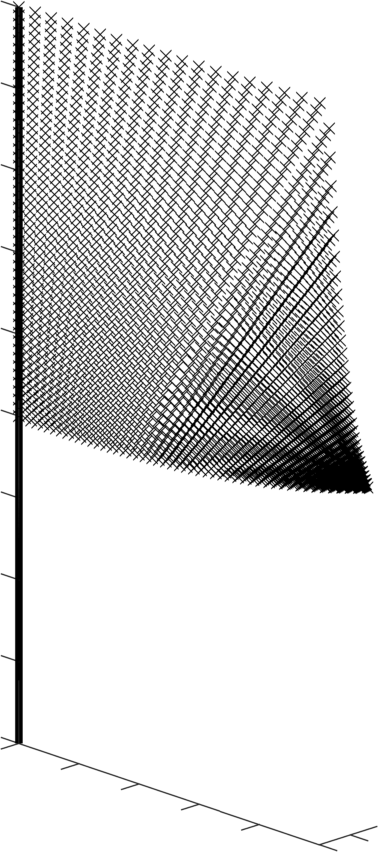}
\qquad
\includegraphics[width=.1\textwidth]{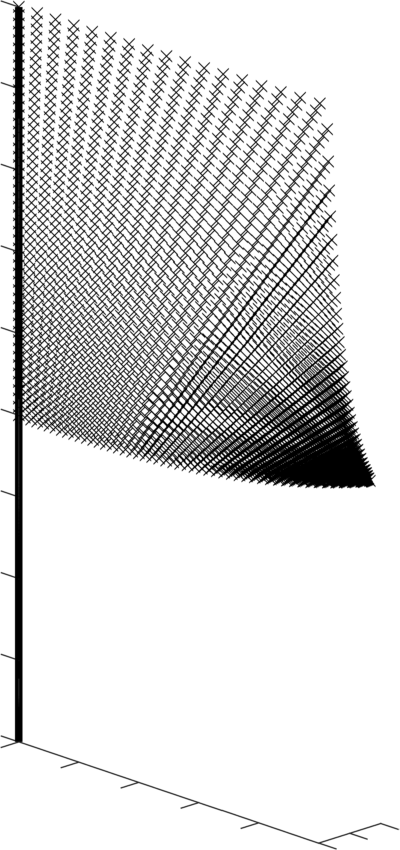}
\caption{A geodesic at times $t=0$, $t=1/3$, $t=2/3$ and $t=1$ in the space of flags. The main part of the deformation is happening in the lower right corner, which is pushed to the back.}
\label{fig:geod:rspace}
\end{figure}

As an example we presentd the solution for the boundary conditions
\begin{equation}
\alpha_0(x)=x, \quad \kappa_0(x)=x,\qquad \alpha_1(x)=1.5x,\quad \kappa_1(x)=\sin(1.5x).
\end{equation}
We discretized the geodesic using 16 timepoints, and 50 points in space, i.e., we solved the minimization problem in 1386 variables, which takes less than one minute on an M1 Macbook  (2021 model). The obtained solution in the space of $\alpha$ and $\kappa$ functions can be seen in Figure~\ref{fig:geod:alphaspace}. In Figure~\ref{fig:geod:rspace} we present the reconstruction of the flag for selected time points of the solution.

\begin{remark}
The difficulty with the numerical approach is partly because the condition \eqref{alphaconstraints2} on $\alpha$ is \emph{not} a convex
condition on the space $\spacey$ of $C^1$ functions $\alpha$. Hence a linear interpolation between two legitimate functions $\alpha_0$ and $\alpha_1$ is unlikely to also satisfy the condition. The other factor leading to difficulty is that the kinetic energy really is smaller if the path in the space
of flags tries to ``cheat'' by making the transformation degenerate. The kinetic energy defined by \eqref{riemannianmetric} involves the area form $du\,dv$, which involves the Jacobian $J(x,y) = 1+y\alpha'(x)$ by formula \eqref{jacobianarea}. The condition on $\alpha$ is equivalent to requiring that $J(x,1) = \lambda(x)$ be positive everywhere, and so $\alpha$ comes closest to violating it when the Jacobian comes close to vanishing. Such a near-degenerate flag configuration will make the kinetic energy smaller for a fixed velocity field. This could be related to vanishing phenomena for $L^2$ Riemannian metrics on spaces of curves or diffeomorphisms, where the same kind of degenerate objects are close in the weak metric to smooth objects and provide shortcuts. See Michor-Mumford~\cite{michor2005vanishing} for an explanation and pictures in those cases.
\end{remark}

\bibliographystyle{abbrv}
\bibliography{references}









\end{document}